\documentclass[fleqn,12pt]{article}
\usepackage{amssymb,amsmath,amsfonts}
\usepackage{color}
\usepackage{graphicx}
\usepackage{latexsym}
\usepackage{amsmath}
\usepackage{amssymb}
\usepackage{graphics}
\usepackage[dvips]{epsfig}
\usepackage[margin=1in]{geometry}

\numberwithin{equation}{section}
\newtheorem{thm}{Theorem}[section]

\newtheorem{lem}[thm]{Lemma}
\newtheorem{prop}[thm]{Proposition}
\newtheorem{rmk}[thm]{Remark}

\newcommand{\ea}{\epsilon}
\newcommand{\ta}{\theta}

\newcommand{\da}{\delta}

\renewcommand{\aa}{\alpha}

\newcommand{\pl}{\partial}

\newcommand{\ga}{\gamma}

\newcommand{\iy}{\infty}

\newcommand{\lt}{\left}
\newcommand{\rt}{\right}
\newcommand{\be}{\begin{equation}}
\newcommand{\bs}{\begin{split}}
\newcommand{\es}{\end{split}}
\newcommand{\ee}{\end{equation}}
\newcommand{\bee}{\begin{equation*}}
\newcommand{\eee}{\end{equation*}}

\newcommand{\ef}{\eqref}
\newcommand{\f}{\frac}

\begin{document}
\begin{center}
\large{ \bf Global Existence of Smooth Solutions and Convergence to Barenblatt Solutions for the Physical Vacuum Free Boundary Problem of Compressible Euler Equations with Damping }
\end{center}
\centerline{ Tao Luo, Huihui Zeng}
\begin{abstract} For the physical vacuum free boundary problem with the sound speed being $C^{{1}/{2}}$-H$\ddot{\rm o}$lder continuous  near vacuum boundaries of the one-dimensional compressible Euler equations with damping,
the  global existence of the smooth solution is proved, which is shown to converge to the Barenblatt self-similar solution for the
the porous media equation with the same total mass when the initial data is a small perturbation of the Barenblatt  solution. The pointwise  convergence with a rate of density, the convergence rate of velocity in supreme norm and the precise expanding rate of the physical vacuum boundaries are also given. The proof is based on a construction of  higher-order  weighted functionals with both space and time
weights capturing the behavior of solutions both near vacuum states and in large time,  an introduction of a new ansatz, higher-order nonlinear energy estimates and elliptic estimates.  \end{abstract}

\section{Introduction}
The aim of this paper is to prove the global existence and time-asymptotic equivalence to the Barenblatt self-similar solutions of smooth solutions for the following physical vacuum free boundary problem for the Euler equations of compressible isentropic flow with damping:
\begin{equation}\label{1.1}
\begin{split}
& \rho_t+(\rho u)_x=0 &  {\rm in}& \ \ {\rm I}(t) : = \left\{(x,t) \lt |   x_-(t) < x < x_+(t), \   t > 0\rt. \right\}, \\
&(\rho u)_t+(p(\rho)+\rho u^2)_x=- \rho u  & {\rm in}& \ \ {\rm I}(t),\\
&\rho>0 &{\rm in }  & \ \ {\rm I} (t),\\
 & \rho=0    &    {\rm on}& \  \ \Gamma(t): = \left\{(x,t) \lt |  x= x_\pm (t) , \   t > 0\rt. \right\},\\
 &   \dot{ \Gamma}(t) =u(\Gamma(t),t ),  & &\\
&(\rho, u)=(\rho_0, u_0) & {\rm on} & \ \  {\rm I}(0):=\left\{(x,t) \lt |   x_-(0) < x < x_+(0), \   t =0\rt. \right\}.
\end{split}
\end{equation}
Here $(x,t)\in \mathbb{R}\times [0,\iy)$,  $\rho $, ${\bf u} $ and $p$ denote, respectively, the space and time variable, density, velocity and pressure; ${\rm I}(t) $, $\Gamma(t)$ and
$\dot{ \Gamma}(t)$ represent, respectively,  the changing domain occupied by the
gas,   moving vacuum boundary and   velocity of $\Gamma(t)$; $-\rho u$ appearing on the right-hand side of $\ef{1.1}_2$ describes the frictional damping.
 We assume that   the pressure satisfies the $\gamma$-law:
$$
p(\rho)= \rho^{\gamma}  \  \ {\rm for} \  \ \gamma>1
$$
(Here the adiabatic constant is set to be unity.) Let $c=\sqrt {p'(\rho)}$ be the sound speed, a vacuum boundary  is called {\it physical} if
$$
0< \left|\frac{\pl c^2}{\pl x}\right|<+\infty $$
in a small neighborhood of the boundary.
In order to capture this physical singularity, the initial density  is supposed to satisfy
   \be\label{156}\begin{split}
\rho_0(x)>0 \ \ {\rm for} \ \ x_-(0)< x<x_+(0), \ \  \rho_0\lt(x_\pm(0)\rt)=0\   \   {\rm and}
 \ \ 0< \lt|\lt(\rho_0^{\ga-1}\rt)_x \lt(x_\pm(0)\rt)\rt| <\iy  .
\end{split} \ee
Let $M\in(0,\iy)$ be the initial total mass, then the conservation law of mass, $\ef{1.1}_1$, gives
$$\int_{x_-(t)}^{x_+(t)}\rho(x,t)dx = \int_{x_-(0)}^{x_+(0)} \rho_0(x)dx =:M  \ \ {\rm for} \ \ t>0.$$

The compressible Euler equations of isentropic flow  with damping is closely related to the
porous media equation (cf. \cite{HL, HMP, HPW, 23}):
\begin{equation}\label{pm}
\rho_t=p(\rho)_{xx},
\end{equation}
when $(\ref{1.1})_2$ is simplified to Darcy's law:
\begin{equation}\label{darcy}
p(\rho)_x=- \rho u.
\end{equation}
For \ef{pm}, basic understanding of the solution   with finite mass is provided by Barenblatt (cf. \cite{ba}), which is given by
\begin{equation}\label{bar}\begin{split}
\bar \rho(x, t)=  (1+  t)^{-\frac{1}{\gamma+1}}\left[A- B(1+ t)^{-\frac{2}{\gamma+1}}x^2\right]^{\frac{1}{\gamma-1}},
\end{split}\end{equation}
where
\begin{equation}\label{mass1}
B=\frac{\ga-1}{2\ga(\ga+1)}  \ \ {\rm and} \ \ A^{\frac{\ga+1}{2(\ga-1)}} = {M} \sqrt{B}  \lt( {\int_{-1}^1 (1-y^2)^{1/(\ga-1)}dy} \rt)^{-1}; \end{equation}
so that the  Barenblatt  self-similar solution defined in ${\rm I}_b(t)$ has the same total mass as that for the solution of \ef{1.1}:
\be\label{massforbaren}
\int_{\bar x_-(t)}^{\bar x_+(t)} \bar\rho(x,t) dx    = M =  \int_{x_-(t)}^{x_+(t)}\rho(x,t)dx  \   \  {\rm for} \ \  t\ge 0,\ee
where
 \begin{equation}\label{IB} {\rm I}_b(t)= \lt\{(x,t)\lt| \ \bar x_-(t)< x< \bar x_+(t),  \  t\ge 0 \rt. \rt\} \ \  {\rm with} \ \   \bar  x_{\pm}(t)=\pm \sqrt{A B^{-1}} (1+ t)^{ {1}/({\gamma+1})}. \end{equation}
The corresponding Barenblatt velocity is defined in ${\rm I}_b(t)$  by
\begin{equation*}\label{}
 \bar u(x,t)=- \frac{p(\bar \rho)_x}{\bar \rho}
 =\f{ x}{(\ga+1)(1+  t)} \  \ {\rm satsifying} \ \ \dot {\bar x}_{\pm}(t)=\bar  u(\bar  x_{\pm}(t),\  t).
\end{equation*}
So, $(\bar\rho,\  \bar u)$  defined in the region ${\rm I}_b (t)$
solves (\ref{pm}) and (\ref{darcy}) and satisfies \ef{massforbaren}.

It is clear that the vacuum boundaries $x=\bar  x_{\pm}(t)$  of Barenblatt's solution are physical. This is the major motivation to study problem \ef{1.1}, the physical vacuum free boundary problem of
compressible Euler equations with damping. To this end,  a class of explicit solutions to problem \ef{1.1} was constructed by Liu in \cite{23}, which are of the following form:
\be\label{liuexplicitsolution}\begin{split} & {\rm I}(t)=\left\{(x,t) \lt |   \  -\sqrt{e(t)/b(t)} < x < \sqrt {e(t)/b(t)}, \ \  t \ge 0\rt. \right\},   \\
& c^2(x, t)={e(t)-b(t)x^2} \ \  {\rm and} \ \   u(x, t)=a(t) x  \ \  {\rm  in } \ \   {\rm I} (t).\end{split}\ee
In \cite{23},  a system of ordinary differential equations for $(e,b,a)(t) $ was derived with $e(t), b(t)>0$ for $t\ge 0$ by substituting \ef {liuexplicitsolution} into $\ef{1.1}_{1, 2}$ and the time-asymptotic equivalence of this explicit solution and  Barenblatt's  solution with the same total mass was shown. Indeed, the Barenblatt solution  of \ef{pm} and \ef{darcy} can be obtained by the same ansatz as \ef{liuexplicitsolution}:
$$
 \bar c^2(x, t)=\bar e(t)-\bar b(t)x^2  \  \  {\rm and } \ \  \bar u(x, t)=\bar a(t)x. $$
 Substituting this into \ef{pm}, \ef{darcy} and \ef{massforbaren} with $\bar {x}_{\pm}(t)=\pm\sqrt{{\bar e(t)}/{\bar b(t)}}$ gives
 $$
   \bar e(t)= \ga A (1+  t)^{-({\gamma-1})/({\gamma+1})},\   \bar b(t)= \ga B(1+ t)^{-1}
  \  \ {\rm and} \  \ \bar a(t)=  {(\gamma+1)^{-1}(1+t)^{-1}},$$
  where $A$ and $B$ are determined by \ef{mass1}. Precisely, it was proved in \cite{23} the following  time-asymptotic equivalence:
  $$
  (a,\ b,\ e)(t)=(\bar a, \ \bar b, \ \bar e)(t)+ O(1)(1+t)^{-1}{\ln (1+t)} \ \ {\rm as}\  \ t\to\infty.$$
A question was raised in \cite{23} whether this equivalence is still true for general solutions to  problem \ef{1.1}.  The purpose of this paper is to prove the global existence of smooth solutions to the physical vacuum free boundary  problem \ef{1.1} for  general initial data which are small perturbations of Barenblatt's   solutions, and the time-asymptotic equivalence of them.  In particular, we obtain the pointwise  convergence with a rate of density which gives the detailed behavior of the density, the convergence rate of velocity in supreme norm and the precise expanding rate of the physical vacuum boundaries.  The results obtained in the present work  also prove the nonlinear asymptotic stability of  Barenblatt's  solutions in the setting of physical vacuum free boundary problems.

The physical vacuum that the sound speed is $C^{ {1}/{2}}$-H$\ddot{\rm o}$lder continuous  across  vacuum boundaries  makes the study of  free boundary problems in compressible fluids challenging and very interesting, even for  the local-in-time existence theory,  because  standard methods of symmetric hyperbolic systems (cf. \cite{17}) do not apply.
Indeed, characteristic speeds of the compressible isentropic Euler equations become singular with infinite spatial derivatives
at  vacuum boundaries which creates much severe  difficulties in analyzing the regularity near boundaries.   The phenomena of physical vacuum arise in several important physical situations naturally besides the above mentioned, for example, the equilibrium and dynamics of boundaries of gaseous stars (cf. \cite{17', LXZ}).
Recently, important progress has been made in the local-in-time well-posedness theory for the one- and three-dimensional compressible
Euler equations (cf. \cite{16, 10, 7, 10', 16'}).
In  the theory and application of nonlinear partial differential equations, it is of fundamental importance to study the global-in-time existence and long time asymptotic behavior  of solutions. However, it  poses a great challenge to extend the local-in-time existence theory to the global one of smooth solutions,  due to the strong degeneracy near vacuum states caused by the singular behavior of  physical vacuum.   The key in analyses is to obtain the
global-in-time regularity of solutions near  vacuum boundaries by establishing the uniform-in-time higher-order estimates, which is nontrivial to achieve due to strong degenerate nonlinear hyperbolic characters.  To the best of our knowledge, the results obtained in this paper are the first ones on the global existence of smooth solutions for the physical vacuum free boundary problems in  inviscid compressible fluids. This is somewhat surprising due to the difficulties mentioned above.
It should be pointed that the $L^p$-convergence of
$L^{\infty}$-weak solutions for the Cauchy problem of the one-dimensional compressible Euler equations with damping  to Barenblatt  solutions of the porous media equations was given in \cite{HMP} with  $p=2$ if $1<\ga\le 2$ and $p=\ga$ if $\ga>2$ and in \cite{HPW} with $p=1$, respectively, using   entropy-type estimates for the solution itself without deriving  estimates for derivatives.
However, the interfaces separating gases and vacuum cannot be traced in the framework of $L^{\infty}$-weak solutions.  The aim of the present work is to understand the behavior and long time dynamics of   physical vacuum boundaries, for which obtaining the global-in-time regularity of solutions is essential.

In order to overcome difficulties in the analysis of obtaining
global-in-time regularities of solutions near vacuum boundaries, we construct higher-order weighted functionals  with both space and time weights, introduce a new ansatz to bypass the obstacle that  Barenblatt solutions do  not solve $\ef{1.1}_2$ exactly which causes errors in large time, and perform higher-order nonlinear energy estimates and elliptic estimates.
In the construction of  higher-order weighted functionals,  the space and time weights are used to capture the behavior of  solutions  near  vacuum states and  ditect the decay of solutions to  Barenblatt solutions, respectively.
The choice of these weights also depends on the behavior both near  vacuum states and in large time of   Barenblatt solutions.
As shown in \cite{16, 10, 7, 10', 16'},
a powerful tool in the study of physical vacuum free boundary problems of nonlinear hyperbolic equations is the weighted energy estimate. It should be remarked that
weighted estimates used in establishing the local-in-time well-posedness theory (cf. \cite{16, 10, 7, 10', 16'}) involve only spatial weights. Yet  weighted estimates only involving spatial weights seem to be limited to proving local existence results. To obtain global-in-time higher-order estimates, we introduce time weights to quantify the large time behavior of solutions.
The choice of time weights may be suggested by looking at the linearized problem to get hints on how the solution decays.
Indeed,  after introducing a new ansatz to correct the error due to the fact that Barenblatt solutions do not solve $\ef{1.1}_2$ exactly, one may decompose the solution of \ef{1.1}
as a sum of the Barenblatt solution, the new ansatz and an error term.
For the linearized problem of the error term around the Barenblatt solution, one may obtain precise time decay rates of weighted norms for various orders of derivatives, while the $L^2$-weighted norm of solution itself is bounded.  However, it requires tremendous efforts to pass from linear analyses to  nonlinear analyses for the nonlinear problem \ef{1.1}.
Our strategy for the nonlinear analysis is using a bootstrap argument for which we identify  an appropriate a priori assumption. This a priori assumption involves not only the $L^{\infty}$-norms of the solution and its first derivatives but also the weighted  $L^{\infty}$-norms of higher-order derivatives with both spatial and temporal weights.
This is one of the new ingredients of the present work  compared with the methods used
either for the local existence theory in \cite{16, 10, 7, 10', 16'} or the nonlinear instability theory in \cite{17'}.
Under this a priori assumption, we  first use elliptic estimates to bound the space-time weighted $L^2$-norms of higher-order derivatives in both normal and tangential-normal directions of vacuum boundaries by the corresponding space-time weighted $L^2$-norms of tangential derivatives. With these bounds, we perform the nonlinear weighted energy estimate by differentiating equations in the tangential direction  to give the uniform in time space-time weighted $L^2$-estimates of various order derivatives in the tangential direction.
It is discovered here that the a priori assumption solely is not enough to close the nonlinear energy estimates, one has to use the bounds obtained in the elliptic estimates also. This gives the uniform in time estimates of the higher-order weighted functional we construct.
The bootstrap argument is closed by verifying the a priori assumption for which we prove the weighted $L^{\infty}$-norms appearing on the a priori assumption can be bounded by  the  the higher-order weighted functional. One of the advantages of our approach is that we can prove the global existence and large time convergence of solutions with the detailed convergence rates simultaneously.  It should be remarked that the convergence rates obtained in this article is the same as those for the linearized problem.

We would like to close this introduction by reviewing some priori results on  vacuum free boundary problems for the compressible Euler equations besides the results mentioned  above.
Some local-in-time well- and ill-posedness
results were obtained in  \cite{JMnew} for the one-dimensional compressible Euler equations for polytropic gases featuring
various  behaviors at the fluid-vacuum interface. In \cite{24}, when the singularity near the vacuum is mild in the sense that $c^\alpha$ is smooth across the interface with  $0 < \alpha \le 1$ for the  sound speed $c$, a  local existence theory  was developed for the one-dimensional Euler equations with damping, based on the adaptation of the theory of symmetric hyperbolic systems  which is not applicable to physical vacuum boundary problems for which only $c^2$, the square of sound speed in stead of $c^{\alpha}$ ( $0 < \alpha \le 1$) , is required to be smooth across the gas-vacuum interface (further development for this can be found in \cite{38}). In \cite{LXZ}, a general uniqueness theorem was proved for three dimensional motions of compressible Euler equations with or without self-gravitation and a new local-in-time well-posedness theory was established for  spherically symmetric motions without imposing the compatibility condition of the first derivative being zero at the center of symmetry. An instability theory of stationary solutions to the physical vacuum free boundary problem for the spherically symmetric  compressible Euler-Poisson equations of gaseous stars as $6/5<\gamma<4/3$ was established in \cite{17'}.  In \cite{zhenlei}, the local-in-time well-posedness of the physical vacuum free boundary problem was investigated for the one-dimensional Euler-Poisson equations, adopting the methods motivated  by those in \cite{10} for the one-dimensional Euler equations.

\section{Reformulation of the problem and main results}

\subsection{Fix the domain and Lagrangian variables}
 We make the initial interval of the Barenblatt solution, $\lt(\bar x_-(0), \ \bar x_+(0)\rt)$,  as the reference interval and define a diffeomorphism
$\eta_0:   \lt(\bar x_-(0), \ \bar x_+(0)\rt)    \to \lt(  x_-(0),  \ x_+(0)\rt)$
by
\bee\label{} \int_{x_-(0)}^{\eta_0(x)}\rho_0(y)dy=\int_{\bar x_-(0)}^{x}\bar\rho_0(y)dy  \ \ {\rm for} \  \ x \in \lt(\bar x_-(0),\ \bar x_+(0)\rt) ,\eee
where $\bar \rho_0(x) : = \bar\rho(x,0) $ is the initial density of the Barenblatt solution. Clearly,
\be\label{2.3}
  \rho_0(\eta_0(x))\eta_{0}'(x) =   \bar \rho_0(x)  \ \ {\rm for} \  \ x \in \lt(\bar x_-(0),\ \bar x_+(0)\rt). \ee
Due to \ef{156}, \ef{bar} and the fact that the total mass of the Barenblatt solution is the same as that of $\rho_0$, \ef{massforbaren}, the diffeomorphism $\eta_0$ is well defined. For simplicity of presentation,  set
$$ \mathcal{I} : = \lt(\bar x_-(0), \ \bar x_+(0)\rt)=\lt(-\sqrt{A /B }, \ \sqrt{A/ B } \rt). $$
 To fix the boundary, we transform system \ef{1.1} into Lagrangian variables. For $x\in \mathcal{I}$, we define the Lagrangian variable $\eta(x, t)$ by
$$
 \eta_t(x, t)= u(\eta(x, t), t) \ \ {\rm for} \  \ t>0 \  \ {\rm and} \ \  \eta(x, 0)=\eta_0(x),
$$
and set the Lagrangian density and velocity  by
\be\label{ldv}f(x, t)=\rho(\eta(x, t), t) \ \  {\rm and} \ \  v(x, t)= u(\eta(x, t), t) .\ee
Then the Lagrangian version of system \ef{1.1} can be written on the reference domain $\mathcal{I} $ as
\be\label{e1-2} \begin{split}
& f_t + f v_x /\eta_x=0    & {\rm in}& \  \    \mathcal{I} \times (0, \iy),\\
& f v_t+\lt( f^\ga  \rt)_x /\eta_x = -fv \ \  &{\rm in}& \  \   \mathcal{I}\times (0, \iy), \\
& (f, v)=(\rho_0(\eta_0), u_0(\eta_0) )  &  {\rm on}& \  \   \mathcal{I} \times \{t=0\}.
\end{split}
\ee
The map $\eta(\cdot, t)$ defined above can be extended to $\bar {\mathcal{I}}= [-\sqrt{A /B }, \ \sqrt{A/ B }  ]$. In the setting, the vacuum free boundaries for problem \ef{1.1} are given by
\be\label{vbs}  x_{\pm}(t)=\eta(\bar x_{\pm}(0), \  t)=\eta\lt(\pm\sqrt{A/ B },t \rt) \ \  {\rm for} \ \  t\ge 0.\ee
It follows from solving $\eqref{e1-2}_1$  and using \ef{2.3}   that
\be\label{ld}
f(x,t)\eta_x(x,t)=\rho_0(\eta_0(x))\eta_{0}'(x)=\bar \rho_0(x), \ x\in \mathcal{I} .
\ee
It should be noticed that we need $\eta_x(x,t)>0$ for $x\in\mathcal{I}$ and $t\ge 0$ to make the Lagrangian transformation sensible, which will be verified in \ef{basic}.
So,  the initial density of the Barenblatt solution, $\bar\rho_0$,  can be regarded as a parameter and system \eqref{e1-2} can be rewritten as
 \be\label{equ}\lt.\begin{split}
&\bar\rho_0 \eta_{tt} + \bar\rho_0 \eta_{t}+ \lt( \bar\rho_0^\ga/\eta_x^\ga  \rt)_x= 0,\ \  &{\rm in}& \  \   \mathcal{I}\times (0, \iy), \\
&(\eta, \eta_t)= \lt( \eta_0,  u_0(\eta_0)\rt),  & {\rm on} & \  \  \mathcal{I} \times\{t=0\}.
\end{split}\rt.\ee

\subsection{Ansatz}
Define the Lagrangian variable $\bar\eta(x, t)$ for the Barenblatt flow in $\bar {\mathcal{I}}$ by
$$
 \bar\eta_t (x, t)= \bar u(\bar\eta(x, t), t)=\f{ \bar\eta(x,t)}{(\ga+1)(1+  t)} \ \ {\rm for} \  \ t>0 \  \ {\rm and} \ \  \bar\eta(x, 0)=x,
$$
so that
\be\label{212}
\bar\eta(x,t)=x(1+ t)^{{1}/({\ga+1})} \ \ {\rm for} \ \  (x,t) \in \bar {\mathcal{I}}\times [0,\iy)
\ee
 and
\bee\label{}\lt.\begin{split}
&\bar \rho_0 \bar\eta_t + \lt(\bar \rho_0^\ga/\bar \eta_x^\ga  \rt)_x= 0 \ \  &{\rm in}& \  \   \mathcal{I}\times (0, \iy) .
\end{split}\rt.\eee
Since $\bar\eta$ does not solve $\ef{equ}_1$ exactly, we introduce a correction $h(t)$ which is the solution of the following initial   value problem of ordinary differential equations:
\be\label{mt}\lt.\begin{split}
  &h_{tt} + h_t   -   {(\bar\eta_x+h)^{-\ga}} /(\ga+1)+  \bar\eta_{xtt}   + \bar\eta_{xt} =0 ,  \ \ t>0, \\
  &h(t=0)=h_t(t=0)=0.
\end{split}\rt.\ee
(Notice that $\bar \eta_x$,  $\bar \eta_{xt}$ and $\bar \eta_{xtt}$ are independent of $x$.)
The new ansatz is then given by
\be\label{ansatz}
\tilde{\eta}(x,t) :=\bar\eta(x,t)+x h(t),
\ee
so that
\be\label{equeta} \begin{split}
&  \bar\rho_0   \tilde\eta_{tt}   +    \bar\rho_0\tilde\eta_t + \lt(\bar \rho_0^\ga/\tilde \eta_x^\ga  \rt)_x= 0  \ \  &{\rm in}& \  \   \mathcal{I}\times (0, \iy).
\end{split}
\ee
It should be noticed that  $\tilde\eta_x$ is independent of $x$. We will prove in the Appendix that $\tilde\eta$ behaves similar to $\bar \eta$, that is, there exist   positive constants $K$ and  $C(n)$ independent of time $t$ such that for all  $t\ge 0$,
\be\label{decay}\begin{split}
&\lt(1 +   t \rt)^{{1}/({\ga+1})} \le \tilde \eta_{x} (t)\le K \lt(1 +   t \rt)^{{1}/({\ga+1})}, \ \  \ \  \tilde\eta_{xt}(t)\ge 0, \\
&\lt|\f{  d^k\tilde \eta_{x}(t)}{dt^k}\rt| \le C(n)\lt(1 +   t \rt)^{\frac{1}{\ga+1}-k},   \ \ k= 1,  \cdots, n.
\end{split}\ee
Moreover,  there exists a certain constant $C$ independent of $t$ such that
\be\label{h} 0\le h(t)\le   (1+t)^{-\frac{\ga}{\ga+1}}\ln(1+t) \ \  {\rm and} \ \
 |h_t|\le (1+t)^{-1-\frac{\ga}{\ga+1}}\ln(1+t), \ \ t\ge 0. \ee
The proof of \ef{h} will also be given in the Appendix.

\subsection{Main results}
Let
$$w(x,t)=\eta(x,t)- \tilde\eta(x,t) .$$
Then, it follows from \ef{equ} and \ef{equeta} that
\be\label{eluerpert}\begin{split}
\bar\rho_0w_{tt}   +
   \bar\rho_0w_t + \lt[ \bar\rho_0^\ga \lt(  (\tilde\eta_x+w_x)^{-\ga} -   {\tilde\eta_x}^{-\ga}  \rt) \rt]_x =0.
\end{split}\ee
Denote
$$ \alpha := 1/(\ga-1), \  \ l:=3+ \min\lt\{m \in  \mathbb{N}: \ \  m> \alpha \rt\}=4 +[\aa] .$$
For $j=0,\cdots, l$ and  $i=0,\cdots, l-j$,  we set
\bee\label{}\begin{split}
\mathcal{ E}_{j}(t) : = & (1+  t)^{2j} \int_\mathcal{I}    \lt[\bar\rho_0\lt(\pl_t^j w\rt)^2 + \bar\rho_0^\ga \lt(\pl_t^j  w_x \rt)^2 + (1+ t) \bar\rho_0   \lt(\pl_t^{j+1} w\rt)^2 \rt] (x, t)  dx  ,  \\
\mathcal{ E}_{j, i}(t): = & (1+  t)^{2j}  \int_\mathcal{I}  \lt[\bar\rho_0^{1+(i-1)(\ga-1) }  \lt(\pl_t^j \pl_x^i w\rt)^2 + \bar\rho_0^{1+(i+1)(\ga-1) } \lt(\pl_t^j  \pl_x^{i+1}w  \rt)^2\rt] (x, t) dx  .
\end{split}\eee
 The higher-order norm is defined by
$$
  \mathcal{E}(t) :=  \sum_{j=0}^l \lt(\mathcal{ E}_{j}(t) + \sum_{i=1}^{l-j} \mathcal{ E}_{j, i}(t)  \rt).
$$
It will be proved in Lemma \ref{lem31} that
$$
  \sup_{x\in \mathcal{I}} \lt\{\sum_{j=0}^3 (1+  t)^{2j} \lt |\pl_t^j w(x, t)\rt |^2   +  \sum_{j=0}^1 (1+  t)^{2j} \lt |\pl_t^j w_x(x, t)\rt|^2\rt\} \le C  \mathcal{E}(t)
$$
for some constant $C$ independent of $t$. So the bound  of $\mathcal{E}(t)$ gives the uniform bound and decay of $w$ and its derivatives. Now, we are ready to state the main result.
{\begin{thm}\label{mainthm} There exists a  constant $\bar \delta >0$ such that if
$\mathcal{E}(0)\le \bar \da,$
then   the  problem \eqref{equ}  admits a global unique smooth solution  in $\mathcal{I}\times[0, \iy)$ satisfying for all $t\ge 0$,
$$\mathcal{E}(t)\le C\mathcal{E}(0) $$
and
\be\label{2}\begin{split}
  & \sup_{x\in \mathcal{I}} \lt\{\sum_{j=0}^3 (1+  t)^{2j} \lt |\pl_t^j w(x, t)\rt |^2   +  \sum_{j=0}^1 (1+  t)^{2j} \lt |\pl_t^j w_x(x, t)\rt |^2\rt\}
   \\
   & \qquad +
  \sup_{x\in \mathcal{I}}\sum_{
  i+j\le l,\   2i+j \ge 4 } (1+  t)^{2j}\lt |  \bar\rho_0^{\f{(\ga-1)(2i+j-3)}{2}}\pl_t^j \pl_x^i w(x, t)\rt|^2  \le C \mathcal{E}(0),
\end{split}\ee
where $C$ is a positive constant  independent of $t$.
\end{thm}
It should be noticed that  the time derivatives involved in the initial higher-order energy norm,
$\mathcal{E}(0)$, can be determined via the equation by the initial data $\rho_0$ and $u_0$  (see \cite{10} for instance).

\begin{rmk} For the linearized equation of \ef{eluerpert} around the Barenblatt solution:
$$
\bar\rho_0w^L_{tt}   +
   \bar\rho_0w^L_t -\ga \lt( \bar\rho_0^\ga \tilde \eta_x^{-(\ga+1)}w^L_x\rt)_x =0,
$$
one may easily show that, for example,  using the weighted energy method and \ef{decay},
$$\sum_{j=0}^{k} \mathcal{ E}_{j}(w^L)(t)\le C \sum_{j=0}^{k} \mathcal{ E}_{j}(w^L)(0),\  \ t\ge 0,$$
for any integer $k\ge 0$. Here  $C>0$ is a constant independent of $t$ and
$$ \mathcal{ E}_{j}(w^L)(t)  := (1+  t)^{2j}  \int_\mathcal{I}  \lt[\bar\rho_0\lt(\pl_t^j w^L\rt)^2 + \bar\rho_0^\ga \lt(\pl_t^j  w^L_x \rt)^2 + (1+t)\bar\rho_0   \lt(\pl_t^{j+1} w^L\rt)^2\rt](x, t)  dx   .$$
In Theorem \ef{mainthm}, we obtain the same decay rates for the nonlinear problem.
\end{rmk}

As a corollary of Theorem \ref{mainthm}, we have the following theorem for solutions to the original vacuum free boundary problem
\ef{1.1}.

{\begin{thm}\label{mainthm1} There exists a  constant $\bar\delta >0$ such that if
 $\mathcal{E}(0)\le \bar\da, $
then  the  problem
\ef{1.1}  admits a global unique smooth solution $\lt(\rho, u,  {\rm I}(t)\rt)$ for $t\in[0,\iy)$  satisfying
 \be\label{1'}\begin{split} \lt|\rho\lt(\eta(x, t),t\rt)-\bar\rho\lt(\bar\eta(x, t), t\rt)\rt|
 \le & C\lt(A-Bx^2\rt)^{\frac{1}{\ga-1}}(1+t)^{-\frac{2}{\ga+1}} \\
 & \times \lt(\sqrt{\mathcal{E}(0)}+(1+t)^{-\frac{\ga}{\ga+1}}\ln(1+t)\rt), \end{split}  \ee
\be\label{2'} \lt|u\lt(\eta(x, t),t\rt)-\bar u\lt(\bar\eta(x, t), t\rt)\rt|\le C(1+t)^{-1}\lt( \sqrt{\mathcal{E}(0)}+(1+t)^{-\frac{\ga}{\ga+1}}\ln(1+t)\rt),  \ee
\be\label{3'} -c_2(1+t)^{\frac{1}{\ga+1}}\le x_-(t)\le  -c_1(1+t)^{\frac{1}{\ga+1}},  \ \ c_1(1+t)^{\frac{1}{\ga+1}}\le x_+(t)\le  c_2(1+t)^{\frac{1}{\ga+1}}, \ee
\be\label{4'} \lt|\frac{d^k x_{\pm}(t)}{dt^k}\rt|\le C(1+t)^{\frac{1}{\ga+1}-k} , \ \   k=1, 2, 3 ,\ee
for all $x \in \mathcal{I}$ and $t\ge 0 $. Here $C$, $c_1$ and $c_2$ are    positive  constants  independent of  $t$.
\end{thm}
The pointwise behavior of the density and the convergence of velocity for the vacuum free boundary problem
\ef{1.1}  to that of the Barenblatt  solution are given by  \ef{1'} and \ef{2'}, respectively. \ef{3'} gives the precise expanding rate
of the vacuum boundaries of the problem \ef{1.1}, which is the same as that for the Barenblatt solution shown in \ef{IB}.
It is also shown in \ef{1'} that
the difference of density   to problem \ef{1.1} and the corresponding Barenblatt density decays at the
rate of $(1+t)^{-{2}/(\ga+1)}$ in $L^\iy$, while the density of  the Barenblatt solution, $\bar\rho$, decays at the
rate of $(1 + t)^{-1/(\ga+1)}$ in $L^\iy$ (see \ef{bar} ).

\section{Proof of Theorem \ref{mainthm}}
The proof  is based on the local existence of smooth solutions  (cf. \cite{10, 16}) and continuation arguments. The uniqueness of the smooth solutions can be obtained as in section 11 of \cite{LXZ}.  In order to prove the global existence
 of  smooth solutions, we need to obtain the uniform-in-time {\it a priori} estimates  on any given time interval $[0, T]$ satisfying $\sup_{t\in [0, T]}\mathcal{E}(t)<\infty$. To this end, we use a bootstrap argument by making the following  {\it a priori} assumption: Let $w$ be a smooth solution to \ef{eluerpert} on $[0 , T]$, there exists a suitably small fixed positive number $\ea_0\in (0,1)$ independent of $t$   such that
\be\label{apriori}\begin{split}
  & \sum_{j=0}^3 (1+  t)^{2j} \lt\|\pl_t^j w(\cdot,t)\rt\|_{L^\iy(\mathcal{I})}^2   +  \sum_{j=0}^1 (1+  t)^{2j} \lt\|\pl_t^j w_x(\cdot,t)\rt\|_{L^\iy(\mathcal{I})}^2
   \\
   & \qquad +
  \sum_{
  i+j\le l,\   2i+j \ge 4 } (1+  t)^{2j}\lt\|  \bar\rho_0^{\f{(\ga-1)(2i+j-3)}{2}}\pl_t^j \pl_x^i w(\cdot,t)\rt\|_{L^\iy(\mathcal{I})}^2 \le \epsilon_0^2 , \ \  t \in [ 0, T].
\end{split}\ee
This in particular implies,  noting  $\ef{decay}$, that for all $\ta\in[0,1]$,
\be\label{basic}
\frac{1}{2}(1+t)^{\frac{1}{\ga+1}}\le \lt(\tilde \eta_x+\ta w_x\rt)(x, t)\le 2K(1+t)^{\frac{1}{\ga+1}}, \  \ (x, t)\in \mathcal{I}\times [0, T],\ee
where $K$ is a positive constant appearing in $\ef{decay}_1$.

Under this {\it a priori} assumption, we prove in section 3.2 the following elliptic estimates:
$$\mathcal{ E}_{j, i}(t) \le C \lt( \widetilde{\mathcal{ E}}_{0}(t) + \sum_{\iota=1}^{i+j}\mathcal{ E}_{\iota}(t)\rt),    \ \ {\rm when}      \ \  j \ge 0,  \ \  i \ge 1, \ \  i+j\le l, $$
where $C$ is a positive constant independent of $t$ and
$$\widetilde{\mathcal{ E}}_{0}(t) =\mathcal{ E}_{0}(t)- \int_\mathcal{I} \bar\rho_0 w^2  (x, t)dx.$$
With the {\it a priori} assumption and elliptic estimates,  we show in section 3.3 the following nonlinear weighted energy estimate: for some  positive constant $C$ independent of $t$,
$$\mathcal{E}_j(t)     \le C \sum_{\iota=0}^j \mathcal{ E}_{\iota}(0), \ \  j=0,1, \cdots, l. $$
Finally, the {\it a priori} assumption \ef{apriori} can be verified in section 3.4 by proving
\begin{align*}
  & \sum_{j=0}^3 (1+  t)^{2j} \lt\|\pl_t^j w(\cdot,t)\rt\|_{L^\iy(\mathcal{I})}^2   +  \sum_{j=0}^1 (1+  t)^{2j} \lt\|\pl_t^j w_x(\cdot,t)\rt\|_{L^\iy(\mathcal{I})}^2
   \\
   & \qquad +
  \sum_{
  i+j\le l,\   2i+j \ge 4 } (1+  t)^{2j}\lt\|  \bar\rho_0^{\f{(\ga-1)(2i+j-3)}{2}} \pl_t^j \pl_x^i w(\cdot,t)\rt\|_{L^\iy(\mathcal{I})}^2  \le C \mathcal{E}(t)
\end{align*}
for some positive constant $C$  independent of $t$.
This closes the whole bootstrap argument for  small initial perturbations  and completes the proof of Theorem \ref{mainthm}.

\subsection{Preliminaries}
In this subsection, we  present some embedding estimates for weighted Sobolev spaces which will be used later  and introduce some notations to simplify the presentation.
Set
$$
 d(x):=dist(x, \partial \mathcal{I})=\min\lt\{x+\sqrt{A/B}, \sqrt{A/B}-x\rt\} , \ \   x \in \mathcal{I}=\lt(-\sqrt{A/B}, \ \sqrt{A/B} \rt).
$$
For any  $a>0$ and nonnegative integer $b$, the  weighted Sobolev space  $H^{a, b}(\mathcal{I})$ is given by
$$ H^{a, b}(\mathcal{I}) := \lt\{   d^{a/2}F\in L^2(\mathcal{I}): \ \  \int_\mathcal{I}     d^a|\pl_x ^k F|^2dx<\infty, \ \  0\le k\le b\rt\}$$
  with the norm
$$ \|F\|^2_{H^{a, b}(\mathcal{I})} := \sum_{k=0}^b \int_\mathcal{I}    d^a|\pl_x^k F|^2dx.$$
Then  for $b\ge  {a}/{2}$, it holds the following {\it embedding of weighted Sobolev spaces} (cf. \cite{18'}):
 $$ H^{a, b}(\mathcal{I})\hookrightarrow H^{b- {a}/{2}}(\mathcal{I})$$
    with the estimate
  \be\label{wsv} \|F\|_{H^{b- {a}/{2}}(\mathcal{I})} \le C(a, b) \|F\|_{H^{a, b}(\mathcal{I})} \ee
for some positive constant $C(a, b)$.

The following general version of the {\it Hardy inequality} whose proof can be found  in \cite{18'} will also be used often in this paper.
 Let $k>1$ be a given real number and $F$ be a function satisfying
$$
\int_0^{\da} x^k\lt(F^2 + F_x^2\rt) dx < \iy,
$$
where $\da$ is a positive constant; then it holds that
$$
\int_0^{\da} x^{k-2} F^2 dx \le C(\delta, k) \int_0^{\da} x^k \lt( F^2 + F_x^2 \rt)  dx,
$$
where $C(\delta, k)$ is a constant depending only on $\da$ and $k$. As a consequence, one has
\bee\label{}
\int_{-\sqrt{A/B}}^0 \lt(x+\sqrt{A/B}\rt)^{k-2} F^2 dx \le C \int_{-\sqrt{A/B}}^0  \lt(x+\sqrt{A/B}\rt)^k \lt( F^2 + F_x^2 \rt)  dx,
\eee
\bee\label{}
\int^{\sqrt{A/B}}_0 \lt(\sqrt{A/B}-x\rt)^{k-2} F^2 dx \le C \int^{\sqrt{A/B}}_0  \lt(\sqrt{A/B}-x\rt)^k \lt( F^2 + F_x^2 \rt)  dx,
\eee
where $C$ is a constant depending on $A$, $B$ and $k$. In particular, it holds that
\be\label{hardybdry}
\int_\mathcal{I}   d(x)^{k-2} F^2 dx \le C \int_\mathcal{I}   d(x)^k \lt( F^2 + F_x^2 \rt)  dx,
\ee
provided that the right-hand side is finite.

\vskip 0.5cm
{\bf Notations:}

1) Throughout the rest of paper,   $C$  will denote a positive constant which  only depend on the parameters of the problem,
$\ga$ and $M$, but does not depend on the data. They are referred as universal and can change
from one inequality to another one. Also we use $C(\beta)$ to denote  a certain positive constant
depending on quantity $\beta$.

2) We will employ the notation $a\lesssim b$ to denote $a\le C b$ and  $a \sim b$ to denote $C^{-1}b\le a\le Cb$,
where  $C$ is the universal constant  as defined
above.

3) In the rest of the paper, we will use the notations
$$ \int=:\int_{\mathcal{I}}  \ , \ \   \|\cdot\|=:\|\cdot\|_{L^2(\mathcal{I})}  \ \    {\rm and} \ \   \|\cdot\|_{L^{\infty}}=:\|\cdot\|_{L^{\infty}(\mathcal{I})}.$$

4) We set
$$
  \varsigma(x):=\bar\rho_0^{\ga-1}(x)=A-Bx^2, \ \  x \in \mathcal{I}.$$
Then $\mathcal{ E}_{j}$ and $\mathcal{ E}_{j, i}$ can be rewritten as
\bee\label{}\begin{split}
&\mathcal{ E}_{j} (t) =   (1+  t)^{2j}  \int  \lt[  \varsigma^{\aa}\lt(\pl_t^j w\rt)^2 +  \varsigma^{\aa+1} \lt(\pl_t^j  w_x \rt)^2    + (1+ t)\varsigma^{\aa}    \lt(\pl_t^{j+1} w\rt)^2\rt]  (x, t)dx      ,\\
&\mathcal{ E}_{j, i}(t) =  (1+  t)^{2j}  \int  \lt[   \varsigma^{\aa+i+1} \lt(\pl_t^j  \pl_x^{i+1}w  \rt)^2  +    \varsigma^{\aa+i-1} \lt(\pl_t^j  \pl_x^{i}w  \rt)^2  \rt] (x, t)dx .
\end{split}\eee
Obviously,  $ \varsigma(x)$ is equivalent to  $d(x)$, that is,
\be\label{varsigma}
B \sqrt{A/B}  d(x)\le \varsigma(x)\le 2 B \sqrt{A/B}  d(x), \ \ x \in \mathcal{I}.
\ee

\subsection{Elliptic estimates }
Denote
\bee\label{}\begin{split}
&\widetilde{\mathcal{ E}}_{0}(t)  :=     \int     \varsigma^{\aa+1} w_x^2(x, t)   dx   + (1+ t)  \int    \varsigma^{\aa}    w_t^2 (x, t)  dx=\mathcal{ E}_{0}(t)- \int \varsigma^{\aa} w^2  (x, t)dx.
\end{split}\eee
 We prove the following elliptic estimates in this subsection.
\begin{prop} Suppose that \ef{apriori} holds for suitably  small positive number $\epsilon_0 \in(0,1) $, then for $0\le t\le T$,
\be\label{ellipticestimate}
\mathcal{ E}_{j, i}(t) \lesssim \widetilde{\mathcal{ E}}_{0}(t) + \sum_{\iota=1}^{i+j}\mathcal{ E}_{\iota}(t),    \ \ {\rm when}      \ \  j \ge 0,  \ \  i \ge 1, \ \  i+j\le l.
\ee
\end{prop}
The proof of this proposition consists of Lemma \ref{lem41} and Lemma \ref{lem43}.

\subsubsection{Lower-order elliptic estimates}
Equation \ef{eluerpert} can be rewritten as
\bee\label{}\begin{split}
\ga {\tilde\eta_x}^{-\ga-1} \lt(\bar\rho_0^\ga w_x\rt)_x=\bar\rho_0w_{tt}   +
   \bar\rho_0w_t + \lt[ \bar\rho_0^\ga \lt(  (\tilde\eta_x+w_x)^{-\ga} -   {\tilde\eta_x}^{-\ga} +\ga {\tilde\eta_x}^{-\ga-1} w_x \rt)   \rt]_x.
\end{split}\eee
Divide the equation above by $\bar\rho_0$ and expand the resulting equation to give
\be\label{007}\begin{split}
\ga  {\tilde\eta_x}^{-\ga-1}\lt[ \varsigma w_{xx} + (\aa+1)  \varsigma_x w_x\rt]
=  &
  w_{tt}   +
    w_t   -
   \ga \varsigma \lt[ (\tilde\eta_x+w_x)^{-\ga-1} -   {\tilde\eta_x}^{-\ga-1}\rt] w_{xx}\\
    &+(1+\aa)    \varsigma_x \lt[  (\tilde\eta_x+w_x)^{-\ga} -   {\tilde\eta_x}^{-\ga} +\ga {\tilde\eta_x}^{-\ga-1} w_x \rt] .
\end{split}\ee
\begin{lem}\label{lem41}
Suppose that \ef{apriori} holds for suitably  small positive number $\epsilon_0 \in(0,1) $. Then,
$$\mathcal{ E}_{0, 1}(t) \lesssim \widetilde{\mathcal{ E}}_0(t) + \mathcal{ E}_1(t), \ 0\le t\le T.$$
\end{lem}
{\bf Proof}. Multiply the equation \ef{007} by $ {\tilde\eta_x}^{\ga+1} \varsigma^{\aa/2}$ and square the spatial $L^2$-norm of the product to obtain
\be\label{n44}\begin{split}
&\lt\|   \varsigma^{1+\f{\aa}{2}} w_{xx} + (\aa+1)   \varsigma^{\f{\aa}{2}} \varsigma_x w_x\rt\|^2
\\
 \lesssim &
 (1+t)^2 \lt( \lt\|   \varsigma^{\f{\aa}{2}} w_{tt}\rt\|^2   +\lt\|
     \varsigma^{\f{\aa}{2}} w_t  \rt\|^2 \rt) +
   {\tilde\eta_x}^{-2} \lt\|   \varsigma^{1+\f{\aa}{2}} w_x w_{xx} \rt\|^2 +
   {\tilde\eta_x}^{-2} \lt\|   \varsigma^{\f{\aa}{2}} \varsigma_x w_x^2 \rt\|^2\\
   \lesssim &  \mathcal{ E}_{1} +  {\tilde\eta_x}^{-2} \|w_x\|_{L^\iy}^2 \lt(\lt\|   \varsigma^{1+\f{\aa}{2}}   w_{xx} \rt\|^2 +
     \lt\|   \varsigma^{\f{\aa}{2}} \varsigma_x w_x  \rt\|^2\rt),
\end{split}\ee
where we have used the Taylor expansion, the smallness of $w_x$ which is the consequence of \ef{apriori}, and \ef{decay} to derive the first inequality and the definition of $\mathcal{ E}_{1}$  the second.
Note that the left-hand side of \ef{n44} can be expanded as
\bee\label{}\begin{split}
&\lt\|   \varsigma^{1+\f{\aa}{2}} w_{xx} + (\aa+1)   \varsigma^{\f{\aa}{2}} \varsigma_x w_x\rt\|^2
\\
 =& \lt\|   \varsigma^{1+\f{\aa}{2}}  w_{xx}\rt\|^2 +(\aa+1)^2 \lt\|   \varsigma^{\f{\aa}{2}} \varsigma_x w_x\rt\|^2
+ (\aa+1) \int    \varsigma^{1+\aa} \varsigma_x \lt(w_x^2\rt)_x    dx \\
= & \lt\|   \varsigma^{1+\f{\aa}{2}}  w_{xx}\rt\|^2 -(\aa+1) \int    \varsigma^{1+\aa} \varsigma_{xx}  w_x^2  dx,
 \end{split}\eee
where the last equality follows from the integration by parts. Thus,
 \be\label{n45}\begin{split}
\lt\|   \varsigma^{1+\f{\aa}{2}}  w_{xx}\rt\|^2 \lesssim   \widetilde{\mathcal{ E}}_{0} +\mathcal{ E}_{1} +  {\tilde\eta_x}^{-2} \|w_x\|_{L^\iy}^2 \lt(\lt\|   \varsigma^{1+\f{\aa}{2}}   w_{xx} \rt\|^2 +
     \lt\|   \varsigma^{\f{\aa}{2}} \varsigma_x w_x  \rt\|^2\rt).
 \end{split}\ee
On the other hand, it follows from \ef{n44} and \ef{n45} that
\bee\label{}\begin{split}
 \lt\| (\aa+1)   \varsigma^{\f{\aa}{2}} \varsigma_x w_x\rt\|^2= & \lt\| \lt[  \varsigma^{1+\f{\aa}{2}} w_{xx} + (\aa+1)   \varsigma^{\f{\aa}{2}} \varsigma_x w_x \rt] -   \varsigma^{1+\f{\aa}{2}} w_{xx} \rt\|^2
\\
\le & 2\lt\|   \varsigma^{1+\f{\aa}{2}} w_{xx} + (\aa+1)   \varsigma^{\f{\aa}{2}} \varsigma_x w_x\rt\|^2
+ 2\lt\|   \varsigma^{1+\f{\aa}{2}}  w_{xx}\rt\|^2\\
\lesssim & \widetilde{\mathcal{ E}}_{0} +\mathcal{ E}_{1}  +  {\tilde\eta_x}^{-2} \|w_x\|_{L^\iy}^2 \lt(\lt\|   \varsigma^{1+\f{\aa}{2}}   w_{xx} \rt\|^2 +
     \lt\|   \varsigma^{\f{\aa}{2}} \varsigma_x w_x  \rt\|^2\rt)
 .
\end{split}\eee
This, together with \ef{n45}, gives
\be\label{similar}\begin{split}
  \lt\|   \varsigma^{1+\f{\aa}{2}}  w_{xx}\rt\|^2
  + \lt\|     \varsigma^{\f{\aa}{2}} \varsigma_x w_x\rt\|^2\lesssim   \widetilde{\mathcal{ E}}_{0} +\mathcal{ E}_{1}  +  {\tilde\eta_x}^{-2} \|w_x\|_{L^\iy}^2 \lt(\lt\|   \varsigma^{1+\f{\aa}{2}}   w_{xx} \rt\|^2 +
     \lt\|   \varsigma^{\f{\aa}{2}} \varsigma_x w_x  \rt\|^2\rt) ,
\end{split}\ee
which implies, with the aid of the smallness of $w_x$ and \ef{decay}, that
\be\label{n46}\begin{split}
  \lt\|   \varsigma^{1+\f{\aa}{2}}  w_{xx}\rt\|^2
  + \lt\|     \varsigma^{\f{\aa}{2}} \varsigma_x w_x\rt\|^2
  \lesssim    \widetilde{\mathcal{ E}}_{0} +\mathcal{ E}_{1}   .
\end{split}\ee
In view of $\lt\|     \varsigma^{{\aa }/{2}}   \varsigma^{1/2} w_x\rt\|^2 \le \widetilde{\mathcal{ E}}_0 $, we then see that $\lt\|     \varsigma^{{\aa }/{2}}   w_x\rt\|^2 \le \widetilde{\mathcal{ E}}_{0} +\mathcal{ E}_{1} $.  Indeed, if we denote
$\mathcal{I}_1=[-\sqrt{A/B}/2,\sqrt{A/B}/2 ]$, then
\be\label{you}\begin{split}
&  \lt\|     \varsigma^{{\aa }/{2}}   w_x\rt\|^2\lesssim \lt\|     \varsigma^{{\aa }/{2}}   w_x\rt\|_{L^2(\mathcal{I}/\mathcal{I}_1)}^2 + \lt\|     \varsigma^{{\aa }/{2}}   w_x\rt\|_{L^2(\mathcal{I}_1)}^2  \\
  \lesssim  & \lt\|     \varsigma^{{\aa }/{2}}  \varsigma_x w_x\rt\|_{L^2(\mathcal{I}/\mathcal{I}_1)}^2 + \lt\|     \varsigma^{{\aa }/{2}}   \varsigma^{1/2}  w_x\rt\|_{L^2(\mathcal{I}_1)}^2
 \le   \lt\|     \varsigma^{{\aa }/{2}}  \varsigma_x w_x\rt\|^2+ \lt\|     \varsigma^{{\aa }/{2}}   \varsigma^{1/2}  w_x\rt\|^2 \lesssim    \widetilde{\mathcal{ E}}_{0} +\mathcal{ E}_{1},
\end{split}\ee
since $\varsigma_x$ and $\varsigma$ have positive lower-order bounds on intervals $\mathcal{I}/\mathcal{I}_1$ and $\mathcal{I}$, respectively.    This completes the proof of Lemma \ref{lem41}. $\Box$

\subsubsection{Higher-order elliptic estimates}
For $i\ge 1$ and $j\ge 0$, $\pl_t^j\pl_x^{i-1}\ef{007}$ yields
\be\label{007t}\begin{split}
 \ga  {\tilde\eta_x}^{-\ga-1}\lt[ \varsigma \pl_t^j \pl_x^{i+1} w + (\aa+i)  \varsigma_x \pl_t^j \pl_x^i w\rt]
=  &
 \pl_t^{j+2} \pl_x^{i-1} w    +
   \pl_t^{j+1} \pl_x^{i-1} w  +Q_1 + Q_2    ,
\end{split}\ee
where
\be\label{q1}\begin{split}
& Q_1 :=-\ga \sum_{\iota=1}^j \lt[\pl_t^\iota \lt({\tilde\eta_x}^{-\ga-1}\rt)\rt]
 \pl_t^{j-\iota}\lt[ \varsigma  \pl_x^{i+1} w + (\aa+i)  \varsigma_x   \pl_x^i w\rt] \\
& -\ga \pl_t^j\lt\{{\tilde\eta_x}^{-\ga-1}
 \lt[\sum_{\iota=2}^{i-1}C_{i-1}^\iota \lt(\pl_x^\iota \varsigma\rt)\lt(\pl_x^{i+1-\iota} w\rt) +(\aa+1)\sum_{\iota=1}^{i-1}C_{i-1}^\iota \lt(\pl_x^{\iota+1} \varsigma\rt)\lt(\pl_x^{i -\iota} w\rt)
\rt] \rt\},
\end{split}\ee
\be\label{q2}\begin{split}
Q_2 : = &   -
   \ga \pl_t^j \pl_x^{i-1} \lt\{ \varsigma \lt[ (\tilde\eta_x+w_x)^{-\ga-1} -   {\tilde\eta_x}^{-\ga-1}\rt] w_{xx}\rt\}\\
    &+(1+\aa)   \pl_t^j \pl_x^{i-1} \lt\{  \varsigma_x \lt[  (\tilde\eta_x+w_x)^{-\ga} -   {\tilde\eta_x}^{-\ga} +\ga {\tilde\eta_x}^{-\ga-1} w_x \rt] \rt\} .
\end{split}\ee
Here and thereafter  $C_m^j$ is used to denote the binomial coefficients for
$0\le j \le m$,
$$C_m^j=\frac{m!}{j!(m-j)!}.$$
In this paper, summations $\sum_{\iota=1}^{i-1}$ and $\sum_{\iota=2}^{i-1}$ should be understood as zero when $i=1$ and $i=1, 2$, respectively. Multiply equation \ef{007t} by ${\tilde\eta_x}^{ \ga+1} \varsigma^{(\aa+i-1)/2}$ and square the spatial $L^2$-norm of the product to give
 \bee\label{}\begin{split}
   \lt\|\varsigma^\frac{\aa+i+1}{2} \pl_t^j \pl_x^{i+1} w + (\aa+i) \varsigma^\frac{\aa+i-1}{2} \varsigma_x \pl_t^j \pl_x^i w \rt\|^2
\lesssim
(1+t)^2 \lt(\lt\| \varsigma^{\frac{\aa+i-1}{2}} \pl_t^{j+2} \pl_x^{i-1} w \rt\|^2
\rt. \\
\lt.   +
 \lt\|\varsigma^{\frac{\aa+i-1}{2}}  \pl_t^{j+1} \pl_x^{i-1} w \rt\|^2 \rt) + (1+t)^2 \lt( \lt\| \varsigma^{\frac{\aa+i-1}{2}} Q_1 \rt\|^2     +
  \lt\|\varsigma^{\frac{\aa+i-1}{2}}  Q_2  \rt\|^2 \rt).
\end{split}\eee
Similar to the derivation of \ef{similar} and \ef{you}, we can then obtain
\be\label{important}\begin{split}
 &(1+t)^{-2j}\mathcal{E}_{j,i}(t) = \lt\|\varsigma^\frac{\aa+i+1}{2} \pl_t^j \pl_x^{i+1} w \rt\|^2 +  \lt\| \varsigma^\frac{\aa+i-1}{2}  \pl_t^j \pl_x^i w \rt\|^2
\lesssim   \lt\| \varsigma^\frac{\aa+i}{2}   \pl_t^j \pl_x^i w \rt\|^2 +
(1+t)^2 \\
&\times \lt(\lt\| \varsigma^{\frac{\aa+i-1}{2}} \pl_t^{j+2} \pl_x^{i-1} w \rt\|^2
 +
 \lt\|\varsigma^{\frac{\aa+i-1}{2}}  \pl_t^{j+1} \pl_x^{i-1} w \rt\|^2 +  \lt\| \varsigma^{\frac{\aa+i-1}{2}} Q_1 \rt\|^2     +
  \lt\|\varsigma^{\frac{\aa+i-1}{2}}  Q_2  \rt\|^2 \rt).
\end{split}\ee
We will use this estimate to prove the following lemma by the mathematical induction.

\begin{lem}\label{lem43} Assume that \ef{apriori} holds for suitably  small positive number $\epsilon_0\in(0,1) $. Then for $ j\ge 0$, $i\ge 1$ and  $2\le i+j\le l$,
\be\label{lem43est}
\mathcal{ E}_{j, i}(t) \lesssim \widetilde{\mathcal{ E}}_{0}(t) + \sum_{\iota=1}^{i+j}\mathcal{ E}_{\iota}(t),   \ \   t\in [0,T].
\ee
\end{lem}
{\bf Proof}. We use the induction for $i+j$ to prove this lemma. As shown in Lemma \ref{lem41} we know that  \ef{lem43est} holds for $i+j=1$. For $1\le k\le l- 1$, we make the induction hypothesis that  \ef{lem43est} holds for all
$ j\ge 0$, $i\ge 1$ and  $i+j\le k$, that is,
\be\label{asssup}
\mathcal{ E}_{j, i}(t) \lesssim \widetilde{\mathcal{ E}}_{0}(t) + \sum_{\iota=1}^{i+j}\mathcal{ E}_{\iota}(t),  \ \  j\ge 0, \ \ i\ge 1, \ \  i+j \le k,
\ee
it then suffices to prove  \ef{lem43est}  for $ j\ge 0$, $i\ge 1$ and  $i+j=k+1$.
(Indeed, there exists an order of $(i,j)$ for the proof. For example,  when $i+j=k+1$ we will bound $\mathcal{ E}_{k+1-\iota, \iota}$ from $\iota=1$ to $k+1$ step by step.)

 We estimate $Q_1$ and $Q_2$ given by \ef{q1} and \ef{q2} as follows.
For $Q_1$, it follows from \ef{decay} that
\bee\label{}\begin{split}
  |Q_1| \lesssim \sum_{\iota=1}^j (1+t)^{-1-\iota}
 \lt( \varsigma \lt|\pl_t^{j-\iota} \pl_x^{i+1} w\rt| +  \lt| \pl_t^{j-\iota}\pl_x^i w \rt|\rt)   +\sum_{\iota=0}^j \sum_{r=1}^{i-1} (1+t)^{-1-\iota} \lt|\pl_t^{j-\iota} \pl_x^{r} w\rt|,
\end{split}\eee
so that
\be\label{building}\begin{split}
  \lt\| \varsigma^{\frac{\aa+i-1}{2}} Q_1 \rt\|^2 \lesssim &\sum_{\iota=1}^j (1+t)^{-2 -2\iota}
\lt(   \lt\| \varsigma^{\frac{\aa+i+1}{2}}  \pl_t^{j-\iota} \pl_x^{i+1} w\rt\|^2 +  \lt\| \varsigma^{\frac{\aa+i-1}{2}} \pl_t^{j-\iota}\pl_x^i w \rt\|^2\rt)   \\
&+\sum_{\iota=0}^j \sum_{r=1}^{i-1} (1+t)^{-2-2\iota} \lt\|\varsigma^{\frac{\aa+i-1}{2}} \pl_t^{j-\iota} \pl_x^{r} w\rt\|^2  \\
\lesssim &  (1+t)^{-2 -2j} \lt(\sum_{\iota=1}^j \mathcal{ E}_{j -\iota, i}+\sum_{\iota=0}^j \sum_{r=1}^{i-1}  \mathcal{ E}_{j -\iota, r} \rt).
\end{split}\ee
For $Q_2$, it follows from \ef{decay} and \ef{apriori} that
\bee\label{}\begin{split}
  |Q_2| \lesssim & \sum_{n=0}^j \sum_{m=0}^{i-1} K_{nm}\lt(\lt|\pl_t^{j-n}\pl_x^{i-1-m} (\varsigma w_{xx})\rt|+ \lt|\pl_t^{j-n}\pl_x^{i-1-m} (\varsigma_x w_{x})\rt|\rt)\\
  \lesssim & \sum_{n=0}^j \sum_{m=0}^{i-1} K_{nm}\lt(\lt|\varsigma \pl_t^{j-n}\pl_x^{i-m+1} w  \rt|+ \lt|\varsigma_x \pl_t^{j-n}\pl_x^{i-m}  w \rt|+ \sum_{r=1}^{i-m-1}  \lt|  \pl_t^{j-n}\pl_x^{r}  w \rt|  \rt) \\
  =: & \sum_{n=0}^j \sum_{m=0}^{i-1} Q_{2nm}.
\end{split}\eee
Here
\bee\label{}\begin{split}
&  K_{00} = \ea_0(1+t)^{-1-\frac{1}{\ga+1}}; \\
&  K_{10}= \ea_0(1+t)^{-2-\frac{1}{\ga+1}}, \ \
  K_{01} = (1+t)^{-1-\frac{1}{\ga+1}}|\pl_x^2 w |;\\
&  K_{20}=\ea_0(1+t)^{-3-\frac{1}{\ga+1}} + (1+t)^{-1-\frac{1}{\ga+1}}\lt|\pl_t^2\pl_x w \rt|,\\
&K_{11}= (1+t)^{-2-\frac{1}{\ga+1}}\lt|\pl_x^2 w  \rt| + (1+t)^{-1-\frac{1}{\ga+1}}\lt|\pl_t\pl_x^2 w  \rt| , \\
&K_{02}=(1+t)^{-1-\frac{1}{\ga+1}}\lt|\pl_x^3 w  \rt| + (1+t)^{-1-\frac{2}{\ga+1}}\lt|\pl_x^2 w  \rt|^2.
\end{split}\eee
We do not list here $K_{nm}$ for $n+m\ge 3$ since we can use the same method to estimate $Q_{2nm}$ for $n+m\ge 3$ as that for $n+m\le 2$. Easily, $Q_{200}$ and $Q_{210}$ can be bounded by
\bee\label{}\begin{split}
\lt\|\varsigma^{\frac{\aa+i-1}{2}}  Q_{200}  \rt\|^2  \lesssim & \ea_0^2(1+t)^{-2 } \lt(\lt\|\varsigma^\frac{\aa+i+1}{2} \pl_t^j \pl_x^{i+1} w \rt\|^2 +  \lt\| \varsigma^\frac{\aa+i-1}{2}  \pl_t^j \pl_x^i w \rt\|^2  \rt.\\
&\lt. + \sum_{r=1}^{i -1}  \lt\|\varsigma^{\frac{\aa+i-1}{2}}  \pl_t^{j }\pl_x^{r}  w \rt\|^2  \rt)
  \lesssim \ea_0^2 (1+t)^{-2-2j} \lt( \mathcal{E}_{j,i} + \sum_{r=1}^{i-1 }\mathcal{E}_{j,r}\rt)  ,
\end{split}\eee
\bee\label{}\begin{split}
\lt\|\varsigma^{\frac{\aa+i-1}{2}}  Q_{2 1 0}  \rt\|^2  \lesssim & \ea_0^2(1+t)^{- 4 } \lt(\lt\|\varsigma^\frac{\aa+i+1}{2} \pl_t^{j-1} \pl_x^{i+1} w \rt\|^2 +  \lt\| \varsigma^\frac{\aa+i-1}{2}  \pl_t^{j-1} \pl_x^i w \rt\|^2  \rt.\\
&\lt. + \sum_{r=1}^{i -1}  \lt\|\varsigma^{\frac{\aa+i-1}{2}}  \pl_t^{j-1 }\pl_x^{r}  w \rt\|^2  \rt)
  \lesssim \ea_0^2(1+t)^{-2-2j}\sum_{r=1}^{i }\mathcal{E}_{j-1,r}.
\end{split}\eee
For $Q_{201}$, we use \ef{apriori} to get $|\varsigma^{1/2} w_{xx}|\lesssim \ea_0$ and then obtain
\bee\label{}\begin{split}
\lt\|\varsigma^{\frac{\aa+i-1}{2}}  Q_{201}  \rt\|^2  \lesssim &  \ea_0^2 (1+t)^{-2} \lt(\lt\|\varsigma^\frac{\aa+i }{2}  \pl_t^{j } \pl_x^{i } w \rt\|^2 +  \lt\| \varsigma^\frac{\aa+i-2}{2}   \pl_t^{j } \pl_x^{i-1} w \rt\|^2  \rt.\\
&\lt. + \sum_{r=1}^{i -2}  \lt\|\varsigma^{\frac{\aa+i-2}{2}}   \pl_t^{j  }\pl_x^{r}  w \rt\|^2  \rt)
  \lesssim \ea_0^2 (1+t)^{-2-2j}\sum_{r=1}^{i-1}\mathcal{E}_{j,r},
\end{split}\eee
It should be noted that $Q_{201}$ appears when $i\ge 2$. Similarly, $Q_{220}$ can be bounded by
\bee\label{}\begin{split}
 \lt\|\varsigma^{\frac{\aa+i-1}{2}}  Q_{2 2 0}  \rt\|^2  \lesssim  &  \ea_0^2(1+t)^{-6} \lt(\lt\|\varsigma^\frac{\aa+i }{2} \pl_t^{j-2} \pl_x^{i+1} w \rt\|^2 +  \lt\| \varsigma^\frac{\aa+i-2}{2}  \pl_t^{j-2} \pl_x^i w \rt\|^2  \rt.\\
& \lt. + \sum_{r=1}^{i -1}  \lt\|\varsigma^{\frac{\aa+i-2}{2}}  \pl_t^{j-2 }\pl_x^{r}  w \rt\|^2  \rt)
\lesssim   \ea_0^2(1+t)^{-2-2j}   \sum_{r=1}^{i+ 1}\mathcal{E}_{j-2,r},
\end{split}\eee
where we have used the Hardy inequality \ef{hardybdry} and the equivalence \ef{varsigma} to derive that
\bee\label{}\begin{split}
  \lt\| \varsigma^\frac{\aa+i-2}{2}  \pl_t^{j-2} \pl_x^i w \rt\|^2
  \lesssim \lt\|\varsigma^\frac{\aa+i}{2} \pl_t^{j-2} \pl_x^{i+1} w \rt\|^2 +  \lt\| \varsigma^\frac{\aa+i }{2}  \pl_t^{j-2} \pl_x^i w \rt\|^2 .
\end{split}\eee
Similar to the estimate  for $ Q_{2 2 0} $, we can obtain
\bee\label{}\begin{split}
\lt\|\varsigma^{\frac{\aa+i-1}{2}}  Q_{211}  \rt\|^2  \lesssim &  \ea_0^2 (1+t)^{-4} \lt(\lt\|\varsigma^\frac{\aa+i-1 }{2}  \pl_t^{j-1 } \pl_x^{i } w \rt\|^2 +  \lt\| \varsigma^\frac{\aa+i-3}{2}   \pl_t^{j-1 } \pl_x^{i-1} w \rt\|^2  \rt.\\
&\lt. + \sum_{r=1}^{i -2}  \lt\|\varsigma^{\frac{\aa+i-3}{2}}   \pl_t^{j -1 }\pl_x^{r}  w \rt\|^2  \rt)
  \lesssim \ea_0^2 (1+t)^{-2-2j}\sum_{r=1}^{i }\mathcal{E}_{j-1,r},
\end{split}\eee
\bee\label{}\begin{split}
\lt\|\varsigma^{\frac{\aa+i-1}{2}}  Q_{202}  \rt\|^2  \lesssim &  \ea_0^2 (1+t)^{-2} \lt(\lt\|\varsigma^\frac{\aa+i-2 }{2}  \pl_t^{j  } \pl_x^{i-1 } w \rt\|^2 +  \lt\| \varsigma^\frac{\aa+i-4}{2}   \pl_t^{j  } \pl_x^{i-2} w \rt\|^2  \rt.\\
&\lt. + \sum_{r=1}^{i -3}  \lt\|\varsigma^{\frac{\aa+i-4}{2}}   \pl_t^{j   }\pl_x^{r}  w \rt\|^2  \rt)
  \lesssim \ea_0^2 (1+t)^{-2-2j}\sum_{r=1}^{i-1}\mathcal{E}_{j,r}.
\end{split}\eee
It should be noted that $Q_{211}$ and $Q_{202}$ appear when $i\ge 2$ and $i\ge 3$, respectively. This ensures the application of the Hardy inequality. Other cases can be done similarly, since the leading term of $K_{nm}$ is $\sum_{\iota=0}^n (1+t)^{-1-\iota-\frac{1}{\ga+1}}\lt|\pl_t^{n-\iota}\pl_x^{m+1} w  \rt|$.
Now, we  conclude that
\bee\label{}\begin{split}
  \lt\| \varsigma^{\frac{\aa+i-1}{2}} Q_2 \rt\|^2 \lesssim &  \ea_0^2(1+t)^{-2-2j}\lt(\mathcal{E}_{j , i} +  \sum_{0\le \iota\le j, \    r \ge 1 , \  \iota+r \le i+j-1 }  \mathcal{ E}_{ \iota, r}  \rt)(t).
\end{split}\eee
Substituting this and \ef{building} into  \ef{important} gives, for suitably small $\ea_0$,
 \be\label{goal}\begin{split}
  \mathcal{E}_{j,i}(t)
\lesssim  &
(1+t)^{2j+2} \lt(\lt\| \varsigma^{\frac{\aa+i-1}{2}} \pl_t^{j+2} \pl_x^{i-1} w \rt\|^2
 +
 \lt\|\varsigma^{\frac{\aa+i-1}{2}}  \pl_t^{j+1} \pl_x^{i-1} w \rt\|^2 \rt) \\
&  + (1+t)^{2j} \lt\| \varsigma^\frac{\aa+i}{2}   \pl_t^j \pl_x^i w \rt\|^2   +   \sum_{0\le \iota\le j, \    r \ge 1 , \  \iota+r \le i+j-1 }   \mathcal{ E}_{ \iota, r}(t)   .
\end{split}\ee
In particular, when $i\ge 3$
 \be\label{goal3}\begin{split}
  \mathcal{E}_{j,i}(t)
\lesssim  & \mathcal{ E}_{ j+2, i-2}  + \mathcal{ E}_{ j+1, i-2} + \mathcal{ E}_{ j, i-1}  + \sum_{0\le \iota\le j, \    r \ge 1 , \  \iota+r \le i+j-1 }   \mathcal{ E}_{ \iota, r}(t)   .
\end{split}\ee

In what follows, we  use \ef{goal} and the induction hypothesis \ef{asssup} to show that \ef{lem43est} holds for $i+j=k+1$.
First, choosing $j=k$ and $i=1$ in \ef{goal}  gives
\bee\label{}\begin{split}
 & \mathcal{E}_{k,1}(t)
\lesssim    \mathcal{E}_{k+1}(t) + \mathcal{E}_{k}(t) +  \sum_{\iota\ge 0, \  r\ge 1, \  \iota+r \le k }  \mathcal{ E}_{ \iota, r}(t)
\end{split}\eee
which, together with \ef{asssup}, implies
\be\label{nek1t}\begin{split}
 & \mathcal{E}_{k,1}(t)
\lesssim    \widetilde{\mathcal{ E}}_{0}(t) + \sum_{\iota=1}^{k+1}\mathcal{ E}_{\iota}(t).
\end{split}\ee
Similarly,
\bee\label{}\begin{split}
 & \mathcal{E}_{k-1,2}(t) \lesssim \mathcal{E}_{k+1}(t) + \mathcal{E}_{k}(t) + \mathcal{E}_{k-1,1}(t) +  \sum_{\iota\ge 0, \  r\ge 1, \  \iota+r \le k }  \mathcal{ E}_{ \iota, r}(t)
\lesssim    \widetilde{\mathcal{ E}}_{0}(t) + \sum_{\iota=1}^{k+1}\mathcal{ E}_{\iota}(t).
\end{split}\eee
For  $\mathcal{E}_{k-2, 3}$, it follows from \ef{goal3}, \ef{nek1t} and  \ef{asssup}  that
\bee\label{}\begin{split}
  \mathcal{E}_{k-2,3}(t)
\lesssim  & \mathcal{ E}_{k , 1}(t)  + \mathcal{ E}_{k -1 , 1}(t) + \mathcal{ E}_{k-2, 2}(t)  + \sum_{\iota\ge 0, \  r\ge 1, \  \iota+r \le k }  \mathcal{ E}_{ \iota, r}(t)  \lesssim    \widetilde{\mathcal{ E}}_{0}(t) + \sum_{\iota=1}^{k+1}\mathcal{ E}_{\iota}(t).
\end{split}\eee
The other cases can be handled similarly. So we have proved \ef{lem43est} when $i+j=k+1$.
This finishes the proof of Lemma \ref{lem43}. $\Box$

\subsection{Nonlinear weighted energy estimates}
In this subsection, we  prove that the weighted energy $\mathcal{E}_j(t)$ can be  bounded by the initial data for $t\in [0,T]$.
\begin{prop} Suppose that \ef{apriori} holds for suitably  small positive number $\epsilon_0 \in(0,1)  $. Then for $t\in [0, T]$,
\be\label{energy}
\mathcal{E}_j(t)     \lesssim \sum_{\iota=0}^j \mathcal{ E}_{\iota}(0), \ \  j=0,1, \cdots, l.
\ee
 \end{prop}
The proof of this proposition consists of Lemma \ref{lem51} and Lemma \ref{lem53}.

\subsubsection{Basic energy estimates}
\begin{lem}\label{lem51} Suppose that \ef{apriori} holds for suitably  small positive number $\epsilon_0 \in(0,1) $. Then,
\be\label{5-5}\begin{split}
  \mathcal{E}_0(t)
 +
 \int_0^t \int \lt[ (1+s) \bar\rho_0w_s^2 +  (1+s)^{-1}  \bar\rho_0^\ga  w_{x}^2 \rt] dxds
  \lesssim   \mathcal{E}_0(0),  \ \ t\in[0,T].
\end{split}\ee
\end{lem}
{\bf Proof}. Multiply \ef{eluerpert} by $w_t$ and integrate the product with respect to the spatial variable  to get
\bee\label{}\begin{split}
\f{d}{dt}\int \f{1}{2}\bar\rho_0w_{t}^2 dx   +
 \int  \bar\rho_0w_t^2 dx - \int \bar\rho_0^\ga \lt[ (\tilde{\eta}_x+w_x)^{-\ga}
   -\tilde{\eta}_x ^{-\ga}  \rt] w_{xt} dx =0.
\end{split}\eee
Note that
\bee\label{}\begin{split}
&\lt[(\tilde{\eta}_x+w_x)^{-\ga} -  \tilde{\eta}_x^{-\ga}  \rt] w_{xt}  \\
 =&\lt[(\tilde{\eta}_x+w_x)^{-\ga}(\tilde{\eta}_x+w_x)_t - (\tilde{\eta}_x+w_x)^{-\ga} \tilde{\eta}_{xt} \rt] -\lt[ \lt(\tilde{\eta}_x^{-\ga}   w_{x}\rt)_t - \lt(\tilde{\eta}_x^{-\ga}\rt)_t w_x\rt] \\
  =& \f{1}{1-\ga}\lt[(\tilde{\eta}_x+w_x)^{1-\ga } - (1-\ga )\tilde{\eta}_x^{-\ga}   w_{x}\rt]_t - \lt[ (\tilde{\eta}_x+w_x)^{-\ga}  +  \ga \tilde{\eta}_x^{-\ga-1}  w_x \rt]\tilde{\eta}_{xt} \\
 =&\f{1}{1-\ga}\lt[(\tilde{\eta}_x+w_x)^{1-\ga } -\tilde{\eta}_x^{1-\ga }-(1-\ga)\tilde{\eta}_x^{-\ga}   w_{x}\rt]_t  - \tilde{\eta}_{xt} \mathfrak{F} ,
\end{split}\eee
where
$$ \mathfrak{F} (x,t):=(\tilde{\eta}_x+w_x)^{-\ga}
   -\tilde{\eta}_x ^{-\ga} +  \ga \tilde{\eta}_x^{-\ga-1}  w_x.
$$
We then have
\be\label{h1}\begin{split}
\f{d}{dt}\int \mathfrak{E}_0(x,t) dx
 +
 \int  \bar\rho_0w_t^2 dx + \int \bar\rho_0^\ga \tilde{\eta}_{xt}  \mathfrak{F} dx=0,
\end{split}\ee
where
 \bee\begin{split}
\mathfrak{E}_0(x,t):=&\f{1}{2}\bar\rho_0w_{t}^2   + \f{1}{\ga-1} \bar\rho_0^\ga \lt[(\tilde{\eta}_x+w_x)^{1-\ga } -\tilde{\eta}_x^{1-\ga }-(1-\ga)\tilde{\eta}_x^{-\ga}   w_{x}\rt]\\
 \sim & \bar\rho_0w_{t}^2 +   \bar\rho_0^\ga \tilde{\eta}_x^{-\ga-1} w_x^2 \sim   \bar\rho_0w_{t}^2 +   {\bar\rho_0^\ga}({1+t})^{-1}  w_x^2 ,
\end{split}\eee
due to the Taylor expansion, the smallness of $|w_x|$ and \ef{decay}. Integrating \ef{h1} gives
\be\label{5-3}\begin{split}
  \int \mathfrak{E}_0(x,t) dx +\int_0^t \int  \bar\rho_0w_s^2 dx ds+  \int_0^t\int\bar\rho_0^\ga \tilde{\eta}_{xs} \mathfrak{F}   dxds
 =  \int \mathfrak{E}_0(x,0) dx.
\end{split}\ee

Multiplying \ef{eluerpert} by $w$ and integrating the product with respect to $x$ and $t$, one has
\bee\label{}\begin{split}
    \lt.\int \bar\rho_0\lt(\f{1}{2} w^2 + ww_t \rt)  dx \rt|_0^t - \int_0^t \int  \bar\rho_0^\ga\lt[ (\tilde{\eta}_{x}+w_x)^{-\ga} -\tilde{\eta}_{x}^{-\ga}   \rt] w_{x}   dx ds =  \int_0^t   \int  \bar\rho_0w_s^2 dx ds .
\end{split}\eee
It follows from the Taylor expansion and the smallness of $|w_x|$ that
$$\lt[ (\tilde{\eta}_{x}+w_x)^{-\ga} -\tilde{\eta}_{x}^{-\ga}   \rt] w_{x} \le -(\ga/2) \tilde{\eta}_{x}^{-\ga-1}  w_{x}^2.$$
We then obtain, using the Cauchy inequality, \ef{5-3} and $\ef{decay} $, that
\be\label{5-4}\begin{split}
  \int \lt(\bar\rho_0 w^2 \rt)(x,t)  dx + \int_0^t \int (1+s)^{-1}  \bar\rho_0^\ga  w_{x}^2   dx ds  \lesssim \int \lt(\bar\rho_0w^2 + \mathfrak{E}_0\rt)(x,0) dx
  =\mathcal{E}_0(0).
\end{split}\ee

Next, we show the time decay of the energy norm.  It follows from \ef{h1} that
\bee\label{}\begin{split}
\f{d}{dt}\lt[(1+t)\int \mathfrak{E}_0(x,t) dx\rt]
 +
 (1+t)\int  \bar\rho_0w_t^2 dx + (1+t)\int \bar\rho_0^\ga \tilde{\eta}_{xt}  \mathfrak{F}_1 dx=\int \mathfrak{E}_0(x,t) dx.
\end{split}\eee
Therefore,
\bee\label{}\begin{split}
 (1+t)\int \mathfrak{E}_0(x,t) dx
 +
 \int_0^t (1+s)\int  \bar\rho_0w_s^2 dxds
\le  \int \mathfrak{E}_0(x,0) dx+\int_0^t \int \mathfrak{E}_0(x,s) dxds \\
\lesssim   \int \mathfrak{E}_0(x,0) dx+\int_0^t \int \lt[  \bar\rho_0w_s^2 +  (1+s)^{-1}  \bar\rho_0^\ga  w_{x}^2   \rt] dxds
\lesssim \mathcal{E}_0(0),
\end{split}\eee
where estimates \ef{5-3} and \ef{5-4} have been used in the  derivation of  the last inequality. This implies
\bee\label{}\begin{split}
 \int \lt[ (1+t)\bar\rho_0w_{t}^2 +   {\bar\rho_0^\ga}   w_x^2 \rt]dx
 +
 \int_0^t (1+s)\int  \bar\rho_0w_s^2 dxds
\lesssim \mathcal{E}_0(0),
\end{split}\eee
which, together with \ef{5-4}, gives \ef {5-5}. This finishes the proof of Lemma \ref{lem51}. $\Box$

\subsubsection{Higher-order energy estimates}
For $k\ge 1$,  $\pl_t^k\ef{eluerpert}$ yields that
\be\label{ptbkt}\begin{split}
\bar\rho_0 \pl_t^{k+2}w  +
   \bar\rho_0\pl_t^{k+1}w -\ga \lt[ {\bar\rho_0^\ga}
     (\tilde{\eta}_{x}+w_x)^{-\ga-1} \pl_t^k w_x + {\bar\rho_0^\ga}J \rt]_x =0,
\end{split}\ee
where
\bee\label{}\begin{split}
    J :=  & \pl_t^{k-1}\lt\{\tilde\eta_{xt}\lt[(\tilde{\eta}_{x}+w_x)^{-\ga-1}-\tilde{\eta}_{x}^{-\ga-1}\rt]
\rt\} \\
& + \lt\{ \pl_t^{k-1}\lt[\lt(\tilde{\eta}_{x}+w_x\rt)^{-\ga-1} w_{xt}\rt]-
    \lt(\tilde{\eta}_{x}+w_x\rt)^{-\ga-1} \pl_t^k w_x \rt\}.
    \end{split}\eee
To obtain the leading terms of $J$,  we single out the terms involving $\pl_t^{k-1} w_x$. To this end, we rewrite $J$ as
\be\label{j1j2}\begin{split}
    J= &  \tilde\eta_{xt}\pl_t^{k-1}\lt[(\tilde{\eta}_{x}+w_x)^{-\ga-1}
    -\tilde{\eta}_{x}^{-\ga-1}\rt]
 +(k-1) \lt[\lt(\tilde{\eta}_{x}+w_x\rt)^{-\ga-1}\rt]_t
    \pl_t^{k-1}w_x \\
 &+ w_{xt} \pl_t^{k-1}\lt[\lt(\tilde{\eta}_{x}+w_x\rt)^{-\ga-1}\rt]
   + \sum_{\iota=1}^{k-1} C_{k-1}^\iota \lt(\pl_t^{\iota}\tilde\eta_{xt}\rt) \pl_t^{k-1-\iota}\lt[(\tilde{\eta}_{x}+w_x)^{-\ga-1}-\tilde{\eta}_{x}^{-\ga-1}\rt]
   \\
 &+ \sum_{\iota=2}^{k-2} C_{k-1}^\iota \lt( \pl_t^{k-\iota}  w_{x} \rt) \pl_t^{\iota}\lt[\lt(\tilde{\eta}_{x}+w_x\rt)^{-\ga-1}\rt]
  \\
  = & k\lt[\tilde{\eta}_{x}+w_x)^{-\ga-1}\rt]_t \pl_t^{k-1}w_x + \tilde{J}
    \end{split}\ee
where
\bee\label{}\begin{split}
 &  \tilde{J} : =      -(\ga+1)   \lt(\tilde\eta_{x }+w_{x }\rt)_t \sum_{\iota=1}^{k-2} C_{k-2}^\iota \lt(\pl_t^{k-1-\iota}w_x\rt) \pl_t^{\iota}\lt[( \tilde{\eta}_{x}+w_x)^{-\ga-2} \rt]  \\
   &  -(\ga+1)  \lt\{ \tilde\eta_{xt}
 \pl_t^{k-2 } \lt\{ \tilde\eta_{xt}\lt[   (\tilde{\eta}_{x}+w_x)^{-\ga-2}  -\tilde{\eta}_{x}^{-\ga-2}
\rt] \rt\} + w_{xt} \pl_t^{k-2 } \lt[(\tilde{\eta}_{x}+w_x)^{-\ga-2} \tilde\eta_{xt}\rt]  \rt\} \\
&+ \sum_{\iota=1}^{k-1} C_{k-1}^\iota \lt(\pl_t^{\iota}\tilde\eta_{xt}\rt) \pl_t^{k-1-\iota}\lt[(\tilde{\eta}_{x}+w_x)^{-\ga-1}-\tilde{\eta}_{x}^{-\ga-1}\rt]
\\
& + \sum_{\iota=2}^{k-2} C_{k-1}^\iota \lt( \pl_t^{k-\iota}  w_{x} \rt) \pl_t^{\iota}\lt[\lt(\tilde{\eta}_{x}+w_x\rt)^{-\ga-1}\rt]
 .
    \end{split}\eee
Here summations $\sum_{\iota=1}^{k-2}$ and $\sum_{\iota=2}^{k-2}$ should be understood as zero when $k=1, 2$ and $k=1, 2, 3$, respectively.
It should be noted  that only the  terms of lower-order derivatives,  $w_x$, $\cdots$, $\pl_t^{k-2}w_x$ are contained in $\tilde J$. In particular, $\tilde J=0$ when $k=1$.

\vskip 1cm
\begin{lem}\label{lem53} Suppose that \ef{apriori} holds for some small positive number $\epsilon_0\in(0,1) $. Then for all   $j= 1,\cdots,l$,
\be\label{new512}\begin{split}
   & \mathcal{E}_j(t)   + \int_0^t   \int \lt[(1+ s)^{2j+1}     \bar\rho_0    \lt(\pl_s^{j+1} w\rt)^2  +  (1+ s)^{2j-1}      \bar\rho_0^\ga    \lt(\pl_s^{j } w_x\rt)^2   \rt] dx ds
  \\
     \lesssim & \sum_{\iota=0}^j \mathcal{ E}_{\iota}(0),     \ \  t\in[0,T].
\end{split}\ee
\end{lem}
{\bf Proof}. We use induction to prove \ef{new512}. As shown in Lemma \ref{lem51}   we know that \ef{new512} holds for $j=0$.
For $1\le k \le l$, we make the induction hypothesis that  \ef{new512} holds for all $j=0, 1, \cdots, k-1$,  i.e.,
\be\label{512}\begin{split}
    & \mathcal{E}_j(t)   + \int_0^t   \int \lt[(1+ s)^{2j+1}     \bar\rho_0    \lt(\pl_s^{j+1} w\rt)^2  +  (1+ s)^{2j-1}      \bar\rho_0^\ga    \lt(\pl_s^{j } w_x\rt)^2   \rt] dx ds\\
     \lesssim & \sum_{\iota=0}^j \mathcal{ E}_{\iota}(0),  \ 0\le j \le k-1.
\end{split}\ee
 It suffices to prove
  \ef{new512} holds for $j=k$ under the induction hypothesis \ef{512} .

{\bf Step 1}. In this step, we prove that
\be\label{516'}\begin{split}
&\f{d}{dt} \int \lt[ \f{1}{2}   \bar\rho_0 \lt(\pl_t^{k+1} w\rt)^2  +  \widetilde{\mathfrak{E}}_k(x,t) \rt] dx   +\int  \bar\rho_0 \lt(\pl_t^{k+1}w\rt)^2 dx \\
  \lesssim &(\da+ \ea_0)  (1+t)^{-2k-2} \mathcal{E}_{k}(t)   + \lt(\da^{-1}+\ea_0\rt) (1+t)^{-2k-2} \lt(\widetilde{\mathcal{E}}_{0} + \sum_{\iota=1}^{k-1}\mathcal{E}_{\iota}\rt)(t)  ,
      \end{split}\ee
for any postive number $\delta>0$ which will be specified later, where
$$
 \widetilde{\mathfrak{E}}_k(x,t) : = \ga {\bar\rho_0^\ga}  \lt[\frac{1}{2} (\tilde{\eta}_{x}+w_x)^{-\ga-1}   \lt(\pl_t^k w_x\rt)^2   +   J  \pl_t^{k } w_x \rt]
$$
satisfying the following estimates:
\be\label{lowbds}\begin{split}
\int  \widetilde{\mathfrak{E}}_k(x,t)dx \ge  C^{-1} (1+t)^{-1}\int {\bar\rho_0^\ga}
 \lt(\pl_t^k w_{x}\rt)^2 dx -C (1+t)^{-2k-1}\lt(\widetilde{\mathcal{E}}_{0} + \sum_{\iota=1}^{k-1}\mathcal{E}_{\iota} \rt)(t),
      \end{split}\ee
\be\label{upbds}\begin{split}
  \int  \widetilde{\mathfrak{E}}_k(x,t)dx \lesssim    (1+t)^{-1}\int {\bar\rho_0^\ga}
 \lt(\pl_t^k w_{x}\rt)^2 dx +  (1+t)^{-2k-1}\lt(\widetilde{\mathcal{E}}_{0} + \sum_{\iota=1}^{k-1}\mathcal{E}_{\iota}\rt)(t).
      \end{split}\ee

We begin with integrating the product of \ef{ptbkt} and $\pl_t^{k+1} w $
 with respect to the spatial variable which gives
\be\label{515}\begin{split}
 \f{d}{dt} \int   \f{1}{2} \bar\rho_0 \lt(\pl_t^{k+1} w\rt)^2     dx   +\int  \bar\rho_0 \lt(\pl_t^{k+1}w\rt)^2 dx  + \f{\ga}{2} \frac{d}{dt}
    \int  {\bar\rho_0^\ga} (\tilde{\eta}_{x}+w_x)^{-\ga-1}   \lt(\pl_t^k w_x\rt)^2 dx
  \\
 = \f{\ga}{2}
    \int  {\bar\rho_0^\ga} \lt[(\tilde{\eta}_{x}+w_x)^{-\ga-1} \rt]_t  \lt(\pl_t^k w_x\rt)^2 dx           - \ga \int  {\bar\rho_0^\ga} J  \pl_t^{k+1} w_x dx.
      \end{split}\ee
We use \ef{j1j2} to estimate the last term on the right-hand side of \ef{515} as follows:
\bee\label{}\begin{split}
& \int  {\bar\rho_0^\ga} J  \pl_t^{k+1} w_x dx =\frac{d}{dt} \int  {\bar\rho_0^\ga} J  \pl_t^{k } w_x dx - \int  {\bar\rho_0^\ga} J_t  \pl_t^{k } w_x dx  \\
 = & \frac{d}{dt} \int  {\bar\rho_0^\ga} J  \pl_t^{k } w_x dx + (\ga+1) k \int {\bar\rho_0^\ga} (\tilde{\eta}_{x}+w_x)^{-\ga-2}\lt(  w_{xt}
 +  \tilde{\eta}_{xt} \rt) \lt(\pl_t^{k  } w_x \rt)^2   dx \\
 & - \int {\bar\rho_0^\ga} \lt\{ k \lt[ (\tilde{\eta}_{x}+w_x)^{-\ga-1}\rt]_{tt}  \pl_t^{k-1} w_x  + \tilde{J}_t  \rt\}
 \pl_t^{k   } w_x      dx.
      \end{split}\eee
It then follows from \ef{515} that
\bee\begin{split}
&\f{d}{dt} \int \lt[ \f{1}{2}   \bar\rho_0 \lt(\pl_t^{k+1} w\rt)^2  +  \widetilde{\mathfrak{E}}_k(x,t) \rt] dx   +\int  \bar\rho_0 \lt(\pl_t^{k+1}w\rt)^2 dx
\\
=&- \ga  (\ga+1)\lt(  k +\frac{1}{2}\rt) \int {\bar\rho_0^\ga} (\tilde{\eta}_{x}+w_x)^{-\ga-2}\lt(  w_{xt}
 +  \tilde{\eta}_{xt} \rt) \lt(\pl_t^{k  } w_x \rt)^2   dx  \\
     & +    \ga \int {\bar\rho_0^\ga} \lt\{ k \lt[ (\tilde{\eta}_{x}+w_x)^{-\ga-1}\rt]_{tt}  \pl_t^{k-1} w_x  + \tilde{J}_t \rt\}
 \pl_t^{k   } w_x      dx .
 \end{split}\eee
This, together with the fact that $\tilde{\eta}_{xt}\ge 0$ and \ef{basic}, implies
\be\label{516}\begin{split}
&\f{d}{dt} \int \lt[ \f{1}{2}   \bar\rho_0 \lt(\pl_t^{k+1} w\rt)^2  +  \widetilde{\mathfrak{E}}_k(x,t) \rt] dx   +\int  \bar\rho_0 \lt(\pl_t^{k+1}w\rt)^2 dx
\\
\le &  -\ga  (\ga+1)\lt(  k +\frac{1}{2}\rt) \int {\bar\rho_0^\ga} (\tilde{\eta}_{x}+w_x)^{-\ga-2}   w_{xt}
   \lt(\pl_t^{k  } w_x \rt)^2   dx \\
& +  k\ga \int {\bar\rho_0^\ga}     \lt[ (\tilde{\eta}_{x}+w_x)^{-\ga-1}\rt]_{tt}  \lt(\pl_t^{k-1} w_x \rt)
 \pl_t^{k   } w_x      dx+ \ga \int {\bar\rho_0^\ga}   \tilde J_t
 \pl_t^{k   } w_x      dx
  .
      \end{split}\ee
For $\widetilde{\mathfrak{E}}_k$,  it follows from  \ef{basic} and the Cauchy-Schwarz inequality that
\be\label{eek1}\begin{split}
  \widetilde{\mathfrak{E}}_k(x,t) \ge   {\bar\rho_0^\ga}\lt[
 C^{-1}(1+t)^{-1}\lt(\pl_t^k w_{x}\rt)^2 -C(1+t)J^2 \rt],
      \end{split}\ee
\be\label{eek2}\begin{split}
  \widetilde{\mathfrak{E}}_k(x,t) \lesssim    {\bar\rho_0^\ga}\lt[
 (1+t)^{-1}\lt(\pl_t^k w_{x}\rt)^2 + (1+t)J^2 \rt].
      \end{split}\ee
It needs to bound $J$, which contains lower-order terms involving $w_x, \cdots, \pl_t^{k-1}w_x$.  Following from \ef{j1j2}, \ef{apriori} and \ef{decay}, one has
\bee\label{}\begin{split}
 J^2 \lesssim & \lt|\lt[ (\tilde{\eta}_{x}+w_x)^{-\ga-1}\rt]_t \rt|^2 \lt( \pl_t^{k-1} w_x \rt)^2   + {\tilde{J}}^2
 \lesssim   (1+t)^{-4}  \lt(\pl_t^{k-1} w_x \rt)^2 + \tilde J^2 \\
 = &(1+t)^{-2k-2}  \lt[(1+t)^{k-1}\pl_t^{k-1} w_x \rt]^2 + \tilde J^2
      \end{split}\eee
which yields that
    \be\label{ggg}\begin{split}
\int  {\bar\rho_0^\ga} (1+t) J^2 dx
\lesssim &(1+t)^{-2k-1} \int {\bar\rho_0^\ga}\lt[(1+t)^{k-1}\pl_t^{k-1} w_x \rt]^2 dx  + (1+t)\int  {\bar\rho_0^\ga}\tilde J^2dx.
      \end{split}\ee
We now show that the integral involving  $\tilde J$ in \ef{ggg}  can   be bounded by $\widetilde{\mathcal{E}}_{0}+ \sum_{\iota=1}^{k-1}\mathcal{E}_{\iota}$ as follows. First, in view of \ef{apriori} and \ef{decay}, we have
\be\label{tildej}\begin{split}
&  | \tilde{J}  | \lesssim     \sum_{\iota=1}^{k-1} (1+t)^{-2-\iota} \lt|\pl_t^{k-1-\iota}w_x\rt|
 +(1+t)^{-1-\frac{1}{ \ga+ 1}}\lt|\pl_t^{2}w_x\rt| \lt|\pl_t^{k-2}w_x\rt| \\
&+(1+t)^{-2-\frac{1}{ \ga+1}}\lt|\pl_t^{2}w_x\rt|\lt|\pl_t^{k-3}w_x\rt| +  (1+t)^{-1-\frac{1}{\ga+1}}\lt|\pl_t^{3}w_x\rt|\lt|\pl_t^{k-3}w_x\rt|
 + {\rm l.o.t.},
\end{split}\ee
where and thereafter  the notation ${\rm l.o.t.}$ is used to represent the lower-order terms involving $\pl_t^{\iota}w_x$ with $\iota=2,\cdots, k-4$.  It should be noticed that the second term on the right-hand side of \ef{tildej} only appears as $k-2\ge 2$, the third term as $k-3\ge 2$ and the fourth as $k-3\ge 3$.
Clearly, we use \ef{apriori} again to obtain
\bee\label{}\begin{split}
   | \tilde{J}  | \lesssim   &  \sum_{\iota=1}^{k-1} (1+t)^{-2-\iota} \lt|\pl_t^{k-1-\iota}w_x\rt|
 +\ea_0 (1+t)^{-3-\frac{1}{ \ga+ 1}} \varsigma^{-\frac{1}{2}}\lt|\pl_t^{k-2}w_x\rt| \\
& + \ea_0 (1+t)^{-4-\frac{1}{\ga+1}}\varsigma^{-1}\lt|\pl_t^{k-3}w_x\rt| + {\rm l.o.t.},
\end{split}\eee
if $k\ge 6$. Similarly, we can bound ${\rm l.o.t.}$ and obtain
\be\label{tildej'}\begin{split}
   | \tilde{J}  | \lesssim   &  \sum_{\iota=1}^{k-1} (1+t)^{-2-\iota} \lt|\pl_t^{k-1-\iota}w_x\rt|
 +\ea_0 \sum_{\iota=2}^{\lt[{k}/{2}\rt]} (1+t)^{-1-\iota-\frac{1}{\ga+1}} \varsigma^{ \frac{1-\iota}{2}}\lt|\pl_t^{k-\iota}w_x\rt| ,
\end{split}\ee
which implies
    \be\label{ggg1}\begin{split}
 (1+t)\int   {\bar\rho_0^\ga}\tilde J^2dx
 \lesssim  &  (1+t)^{-2k-1}\lt(\widetilde{\mathcal{E}}_0(t) + \sum_{\iota=1}^{k-2}   \mathcal{E}_{k-1-\iota}(t)\rt) \\
& + \ea_0^2 \sum_{\iota=2}^{\lt[{k}/{2}\rt]}  (1+t)^{-1-2\iota-\frac{2}{\ga+1}} \int \bar\rho_0^\ga  \varsigma^{ {1-\iota} }\lt|\pl_t^{k-\iota}w_x\rt|^2 dx .
      \end{split}\ee
In view of the Hardy inequality \ef{hardybdry} and the equivalence of $d$ and $\varsigma$ , \ef{varsigma}, we see that for $\iota=2,\cdots,\lt[{k}/{2}\rt]$,
    \bee\label{}\begin{split}
 & \int \bar\rho_0^\ga  \varsigma^{ {1-\iota} }\lt|\pl_t^{k-\iota}w_x\rt|^2 dx
  =\int   \varsigma^{ {\aa + 2 -\iota} }\lt|\pl_t^{k-\iota}w_x\rt|^2 dx  \\
  \lesssim &  \int   \varsigma^{  \aa + 2 -\iota +2  }\lt(\lt|\pl_t^{k-\iota}w_x\rt|^2
  + \lt|\pl_t^{k-\iota}w_{xx}\rt|^2\rt)dx   \\
  \lesssim & \sum_{i=0}^{\iota-1}  \int   \varsigma^{  \aa +  \iota  } \lt|\pl_t^{k-\iota}\pl_x^{i}w_x\rt|^2 dx \lesssim (1+t)^{2\iota-2k}\sum_{i=1}^{\iota-1} \mathcal{E}_{k-\iota, i}.
      \end{split}\eee
Since  $\aa + 2 -\iota \ge \aa - \lt[\lt[\aa\rt]/{2}\rt] \ge 0$ for  $k\le l$, which  ensures the application of the  Hardy inequality.
It then yields from \ef{ggg}, \ef{ggg1} and the elliptic estimates \ef{ellipticestimate} that
    \be\label{jjj}\begin{split}
 (1+t)\int  {\bar\rho_0^\ga} J^2dx
 \lesssim & (1+t)^{-2k-1}\lt(\widetilde{\mathcal{E}}_{0}(t)+ \sum_{\iota=1}^{k-1}\mathcal{E}_{\iota}(t)\rt).
      \end{split}\ee
This, together with \ef{eek1} and \ef{eek2}, proves \ef{lowbds} and \ef{upbds}.

In what follows, we estimate the terms on the right-hand side of \ef{516} to prove \ef{516'}.  It follows from \ef{apriori} and \ef{decay} that
\bee\label{}\begin{split}
\lt|\int {\bar\rho_0^\ga} (\tilde{\eta}_{x}+w_x)^{-\ga-2}   w_{xt}
   \lt(\pl_t^{k  } w_x \rt)^2   dx \rt| \lesssim
   \ea_0 (1+t)^{-2-\frac{1}{\ga+1}}\int {\bar\rho_0^\ga}
   \lt(\pl_t^{k  } w_x \rt)^2   dx
      \end{split}\eee
and
\bee\label{}\begin{split}
&\lt| \int {\bar\rho_0^\ga}   \lt[ (\tilde{\eta}_{x}+w_x)^{-\ga-1}\rt]_{tt} \lt( \pl_t^{k-1} w_x \rt)   \pl_t^k w_x      dx \rt|\\
 \lesssim &    \int {\bar\rho_0^\ga}   \lt[(1+t)^{-3}  +(1+t)^{-1-\frac{1}{\ga+1}}|w_{xtt}|   \rt]\lt|\pl_t^{k-1} w_x \rt| \lt|  \pl_t^k w_x \rt|      dx \\
 \lesssim  & (1+t)^{-3} \int   {\bar\rho_0^\ga}   \lt| \pl_t^{k-1} w_x \rt| \lt|    \pl_t^k w_x \rt|      dx    + \ea_0 (1+t)^{-3-\frac{1}{\ga+1}} \int \bar\rho_0^{(\ga+1)/2}    \lt| \pl_t^{k-1} w_x \rt| \lt|    \pl_t^k w_x \rt|      dx \\ \lesssim  & \lt(\da +\ea_0\rt) (1+t)^{-2} \int   {\bar\rho_0^\ga}   \lt(\pl_t^k w_x\rt)^2   dx  + \da^{-1} (1+t)^{-4} \int   {\bar\rho_0^\ga}   \lt(\pl_t^{k-1} w_x\rt)^2   dx \\
 & + \ea_0 (1+t)^{-4-\frac{2}{\ga+1}} \int   {\bar\rho_0 }   \lt(\pl_t^{k-1} w_x\rt)^2   dx,
      \end{split}\eee
for any $\da>0$. With the help of elliptic estimates \ef{ellipticestimate}, we have
\bee\label{}\begin{split}
\int   {\bar\rho_0 }   \lt(\pl_t^{k-1} w_x\rt)^2   dx
\le & (1+t)^{-2(k-1)}\mathcal{E}_{k-1,1}(t)
\lesssim (1+t)^{-2(k-1)}  \lt(\widetilde{\mathcal{E}}_{0} + \sum_{\iota=1}^{k }\mathcal{E}_{\iota}\rt)(t).
      \end{split}\eee
Therefore, we obtain \ef{516'}  from \ef{516};
since the last integral  on the right-hand side of \ef{516} can be bounded by
\be\label{uu}\begin{split}
 \int {\bar\rho_0^\ga} \lt|  \tilde J_t
 \pl_t^{k   } w_x \rt|     dx    \lesssim & (\da+ \ea_0)  (1+t)^{-2k-2} \mathcal{E}_{k}(t)   \\& + \lt(\da^{-1}+\ea_0\rt)  (1+t)^{-2k-2} \lt(\widetilde{\mathcal{E}}_{0} + \sum_{\iota=1}^{k-1}\mathcal{E}_{\iota}\rt)(t).
      \end{split}\ee

It remains to prove \ef{uu}. Indeed, we   obtain in a similar way to derving \ef{tildej'} that
\bee\label{}\begin{split}
   | \tilde{J}_t  | \lesssim  &  \sum_{\iota=2}^{k} (1+t)^{-2-\iota} \lt|\pl_t^{k-\iota}w_x\rt|  + \lt[(1+t)^{-1-\frac{1}{ \ga+1}}\lt|\pl_t^{3}w_x \rt|
+(1+t)^{-2-\frac{1}{ \ga+1}}\lt|\pl_t^{2}w_x\rt|\rt]\\
&\times\lt|\pl_t^{k-2}w_x\rt|
 + \lt[(1+t)^{-1-\frac{1}{ \ga+1}}\lt|\pl_t^{4}w_x\rt|
+(1+t)^{-2-\frac{1}{ \ga+1}}\lt|\pl_t^{3}w_x\rt|
\rt.\\
&\lt. + (1+t)^{-3-\frac{1}{ \ga+1}}\lt|\pl_t^{2}w_x\rt|
+(1+t)^{-1-\frac{2}{ \ga+1}}\lt|\pl_t^{2}w_x\rt|^2 \rt]\lt|\pl_t^{k-3}w_x\rt| + {\rm l.o.t.}
\end{split}\eee
and
 \bee\label{}\begin{split}
   | \tilde{J}_t  | \lesssim   &  \sum_{\iota=2}^{k} (1+t)^{-2-\iota} \lt|\pl_t^{k-\iota}w_x\rt|
 +\ea_0  \sum_{\iota=2}^{\lt[({k-1})/{2}\rt] } (1+t)^{-2-\iota-\frac{1}{\ga+1}}   \varsigma^{- \frac{ \iota}{2}} \lt|\pl_t^{k-\iota}w_x\rt| ,
\end{split}\eee
which implies
 \bee\label{}\begin{split}
  \int {\bar\rho_0^\ga} & \lt|  \tilde J_t
 \pl_t^{k   } w_x \rt|     dx    \lesssim   \sum_{\iota=2}^{k} (1+t)^{-2-\iota} \int   {\bar\rho_0^\ga} \lt|\pl_t^{k-\iota}w_x\rt| \lt|
 \pl_t^{k   } w_x \rt| dx \\
& + \ea_0 \sum_{\iota=2}^{\lt[({k-1})/{2}\rt]} (1+t)^{-2-\iota-\frac{1}{\ga+1}} \int \bar\rho^\ga \varsigma^{- \frac{ \iota}{2}}\lt|\pl_t^{k-\iota}w_x\rt| \lt|
 \pl_t^{k   } w_x \rt| dx =:P_1+P_2.
      \end{split}\eee
Easily,  it follows from the Cauchy-Schwarz inequality that for any $\da>0$,
\bee\label{}\begin{split}
 P_1   \lesssim & \da(1+t)^{-2 }\int {\bar\rho_0^\ga} \lt|\pl_t^{k }w_x\rt|^2 dx
 + \da^{-1}\sum_{\iota=2}^{k}(1+t)^{-2-2\iota}\int {\bar\rho_0^\ga} \lt|\pl_t^{k-\iota}w_x\rt|^2 dx \\
 \le & \da (1+t)^{-2k-2} \mathcal{E}_{k}(t)   +  \da^{-1} (1+t)^{-2k-2}\lt( \widetilde{\mathcal{E}}_{0} + \sum_{\iota=1}^{k-2} \mathcal{E}_{\iota}\rt)(t)      \end{split}\eee
      and
\bee\label{}\begin{split}
 P_2   \lesssim  & \ea_0 (1+t)^{-2 }\int {\bar\rho_0^\ga} \lt|\pl_t^{k }w_x\rt|^2 dx
 + \ea_0 \sum_{\iota=2}^{\lt[({k-1})/{2}\rt]} (1+t)^{-2-2\iota }  \int {\bar\rho_0^\ga} \varsigma^{-\iota} \lt|\pl_t^{k-\iota}w_x\rt|^2 dx \\
 \lesssim & \ea_0 (1+t)^{-2k-2}\lt(\widetilde{\mathcal{E}}_{0} + \sum_{\iota=1}^{k }\mathcal{E}_{\iota}\rt)(t).
      \end{split}\eee
Here the last inequality comes from the following estimate:
for $\iota=2,\cdots,\lt[({k-1})/{2}\rt]$,
    \bee\label{}\begin{split}
   \int \bar\rho^\ga  \varsigma^{ { -\iota} }\lt|\pl_t^{k-\iota}w_x\rt|^2 dx
  = & \int   \varsigma^{ {\aa + 1 -\iota} }\lt|\pl_t^{k-\iota}w_x\rt|^2 dx
  \lesssim   \sum_{i=0}^{\iota }  \int   \varsigma^{  \aa  +1 +  \iota  } \lt|\pl_t^{k-\iota}\pl_x^{i}w_x\rt|^2 dx
\\
  &\lesssim (1+t)^{2\iota-2k}\sum_{i=1}^{\iota } \mathcal{E}_{k-\iota, i} \lesssim (1+t)^{2\iota-2k} \lt(\widetilde{\mathcal{E}}_{0} + \sum_{\iota=1}^{k }\mathcal{E}_{\iota}\rt)(t),
      \end{split}\eee
which is deduced from  the Hardy inequality \ef{hardybdry}, the equivalence   \ef{varsigma} and the elliptic estimate \ef{ellipticestimate}  where we note that  $\aa + 1 -\iota \ge \aa - \lt[([\aa]+1)/2\rt] \ge 0 $ for  $k\le l$ so that   the  Hardy inequality can be applied.
Now, we finish the  proof of \ef{uu} and obtain \ef{516'}.

\vskip 0.5cm
\noindent{\bf Step 2}. In this step, we prove that
\be\label{521}\begin{split}
&\f{d}{dt} \int  {\mathfrak{E}}_k(x,t)   dx   +\int  \bar\rho_0 \lt(\pl_t^{k+1}w\rt)^2 dx   +
  (1+t)^{-1}  \int  {\bar\rho_0^\ga}   \lt(\pl_t^k w_x\rt)^2 dx \\
  \lesssim &  (1+t)^{-2k-1} \lt(\widetilde{\mathcal{E}}_{0} + \sum_{\iota=1}^{k-1}\mathcal{E}_{\iota}\rt)(t)
  + (1+t)^{-2} \int \bar\rho_0 \lt(\pl_t^{k }w\rt)^2 dx  \\
  \lesssim & (1+t)^{-2k-1} \lt(\widetilde{\mathcal{E}}_{0} + \sum_{\iota=1}^{k-1}\mathcal{E}_{\iota}\rt)(t),
      \end{split}\ee
 where
\bee\label{}\begin{split}
  {\mathfrak{E}}_k(x,t): = \bar\rho_0 \lt[ \f{1}{2} \lt(\pl_t^k w\rt)^2  + \lt(\pl_t^{k+1}w\rt) \pl_t^k w   +    \lt(\pl_t^{k+1} w\rt)^2 \rt]+  2\widetilde{\mathfrak{E}}_k(x,t).
      \end{split}\eee
We start with  integrating the product of \ef{ptbkt} and $\pl_t^k w $ with respect  to $x$ to yield
\be\label{514}\begin{split}
&\f{d}{dt} \int \bar\rho_0\lt(  \lt(\pl_t^{k+1}w\rt) \pl_t^k w + \f{1}{2} \lt(\pl_t^k w\rt)^2     \rt) dx  +\ga
    \int  {\bar\rho_0^\ga} (\tilde{\eta}_{x}+w_x)^{-\ga-1} \lt(\pl_t^k w_x\rt)^2 dx \\
    = & \int  \bar\rho_0 \lt(\pl_t^{k+1}w\rt)^2 dx
   -
    \ga \int  {\bar\rho_0^\ga} J  \pl_t^k w_x dx.
\end{split}\ee
It follows from  $\ef{514}+2\times\ef{516'}$ that
\be\label{519}\begin{split}
&\f{d}{dt} \int  {\mathfrak{E}}_k(x,t)   dx   +\int  \bar\rho_0 \lt(\pl_t^{k+1}w\rt)^2 dx   +\ga
    \int  {\bar\rho_0^\ga} (\tilde{\eta}_{x}+w_x)^{-\ga-1} \lt(\pl_t^k w_x\rt)^2 dx \\
  \lesssim &
    \lt|\int  {\bar\rho_0^\ga} J  \pl_t^k w_x dx \rt| + (\da+ \ea_0)  (1+t)^{-2k-2} \mathcal{E}_{k}(t)   \\
    &+ \lt(\da^{-1}+\ea_0\rt) (1+t)^{-2k-2} \lt(\widetilde{\mathcal{E}}_{0} + \sum_{\iota=1}^{k-1}\mathcal{E}_{\iota}\rt)(t),
      \end{split}\ee
For the first term on the right-hand side of \ef{519}, we have, with the aid of the
Cauchy-Schwarz inequality and \ef{jjj}, that
\bee\label{}\begin{split}
 \lt|\int  {\bar\rho_0^\ga} J  \pl_t^k w_x dx \rt|
 \lesssim    \da (1+t)^{-1} \int   {\bar\rho_0^\ga}   \lt(\pl_t^k w_x\rt)^2   dx  +  \da^{-1} (1+t)\int  {\bar\rho_0^\ga} J^2dx \\
 \lesssim   \da (1+t)^{-1} \int   {\bar\rho_0^\ga}   \lt(\pl_t^k w_x\rt)^2   dx  +  \da^{-1}
    (1+t)^{-2k-1} \lt(\widetilde{\mathcal{E}}_{0} + \sum_{\iota=1}^{k-1}\mathcal{E}_{\iota}\rt)(t).  \end{split}\eee
This finishes the proof of \ef{521}, by using \ef{basic}, noting the smallness of $\ea_0$ and choosing $\da$ suitably small. Moreover,  we deduce from \ef{lowbds} and \ef{upbds}  that
\be\label{fl1}\begin{split}
  \int {\mathfrak{E}}_k(x,t)dx \ge &C^{-1} \int \bar\rho_0 \lt[  \lt(\pl_t^k w\rt)^2    +    \lt(\pl_t^{k+1} w\rt)^2 \rt]dx+  C^{-1} (1+t)^{-1}\int {\bar\rho_0^\ga}
 \lt(\pl_t^k w_{x}\rt)^2 dx \\
  &-C (1+t)^{-2k-1} \lt(\widetilde{\mathcal{E}}_{0} + \sum_{\iota=1}^{k-1}\mathcal{E}_{\iota}\rt)(t),
      \end{split}\ee
\be\label{fl2}\begin{split}
  \int {\mathfrak{E}}_k(x,t)dx \lesssim & \int \bar\rho_0 \lt[  \lt(\pl_t^k w\rt)^2    +    \lt(\pl_t^{k+1} w\rt)^2 \rt]dx+    (1+t)^{-1}\int {\bar\rho_0^\ga}
 \lt(\pl_t^k w_{x}\rt)^2 dx \\
 &+  (1+t)^{-2k-1} \lt(\widetilde{\mathcal{E}}_{0} + \sum_{\iota=1}^{k-1}\mathcal{E}_{\iota}\rt)(t).
      \end{split}\ee

\vskip 0.5cm
{\bf Step 3}. We show the time decay of the norm in this step.
We integrate  \ef{521} and use  the induction hypothesis \ef{512} to show,  noting  \ef{fl1} and \ef{fl2}, that
\bee\label{}\begin{split}
 & \int \lt[\bar\rho_0\lt(  \lt(\pl_t^k w\rt)^2    +    \lt(\pl_t^{k+1} w\rt)^2  \rt)
+
 {\bar\rho_0^\ga}(1+t)^{-1}\lt(\pl_t^k w_x\rt)^2 \rt](x,t) dx  \\
 &  +\int_0^t  \int  \bar\rho_0 \lt(\pl_s^{k+1}w\rt)^2 dx ds+ \int_0^t (1+s)^{-1} \int  {\bar\rho_0^\ga}  \lt(\pl_s^{k}w_x\rt)^2 dxds
 \\ \lesssim  & \sum_{\iota=0}^{k} \mathcal{E}_{\iota}(0) +  \int_0^t    (1+s)^{-2k-1} \lt(\widetilde{\mathcal{E}}_{0} + \sum_{\iota=1}^{k-1}\mathcal{E}_{\iota}\rt)(s)ds
   \lesssim \sum_{\iota=0}^{k} \mathcal{E}_{\iota}(0).
\end{split}\eee
Here the following estimate has been used,
\bee\label{}\begin{split}
&\int_0^t    (1+s)^{ -1} \lt(\widetilde{\mathcal{E}}_{0} + \sum_{\iota=1}^{k-1}\mathcal{E}_{\iota}\rt)(s)ds =
  \int_0^t  \int   \lt[ (1+ s)^{ -1} \bar\rho_0^\ga   w_x^2 +  \bar\rho_0 w_s^2 \rt] dx ds
\\
&+\sum_{\iota=1}^{k-1}\int_0^t (1+ s)^{2(\iota-1)+1}  \int  \bar\rho_0  \lt(\pl_s^\iota w\rt)^2  dxds  + \sum_{\iota=1}^{k-1} \int_0^t (1+ s)^{2\iota-1}  \int   \bar\rho_0^\ga \lt(\pl_s^\iota  w_x \rt)^2  dx ds
\\
&+\sum_{\iota=1}^{k-1}\int_0^t  (1+ s)^ {2\iota}  \int  \bar\rho_0    \lt(\pl_s^{\iota+1} w\rt)^2   dx
   \lesssim \sum_{\iota=0}^{k} \mathcal{E}_{\iota}(0).
\end{split}\eee
Multiplying \ef{521} by $(1+t)^{p}$  and integrating the product with respect to the temporal variable from $p=1$ to $p=2k$ step by step, we can get
\be\label{ai1}\begin{split}
 & (1+t)^{2k} \int \lt[\bar\rho_0\lt(  \lt(\pl_t^k w\rt)^2    +    \lt(\pl_t^{k+1} w\rt)^2  \rt)
+
 {\bar\rho_0^\ga}(1+t)^{-1}\lt(\pl_t^k w_x\rt)^2 \rt](x,t) dx  \\
  & +\int_0^t (1+s)^{2k} \int  \bar\rho_0 \lt(\pl_s^{k+1}w\rt)^2 dx ds + \int_0^t (1+s)^{2k-1} \int  {\bar\rho_0^\ga}  \lt(\pl_s^{k}w_x\rt)^2 dxds \\
  \lesssim & \sum_{\iota=0}^{k} \mathcal{E}_{\iota}(0) +  \int_0^t    (1+s)^{-1} \lt(\widetilde{\mathcal{E}}_{0} + \sum_{\iota=1}^{k-1}\mathcal{E}_{\iota}\rt)(s)ds
   \lesssim \sum_{\iota=0}^{k} \mathcal{E}_{\iota}(0) .
\end{split}\ee
With this estimate at hand, we  multiply \ef{516'} by $(1+t)^{2k+1}$ and  integrate  the product with respect to the temporal variable to get, in view of \ef{ai1} and \ef{512},
\be\label{ai2}\begin{split}
 & (1+t)^{2k+1} \int \lt[\bar\rho_0     \lt(\pl_t^{k+1} w\rt)^2
+
 {\bar\rho_0^\ga}(1+t)^{-1}\lt(\pl_t^k w_x\rt)^2 \rt] dx  \\
 & +\int_0^t   \int  (1+ s)^{2k+1}     \bar\rho_0    \lt(\pl_s^{k+1} w\rt)^2    dx ds\\
 \lesssim &    \sum_{\iota=0}^{k} \mathcal{E}_{\iota}(0) +    \int_0^t (1+s)^{-1} \int \mathcal{E}_{k}(s)ds + \int_0^t    (1+s)^{-1} \lt(\widetilde{\mathcal{E}}_{0} + \sum_{\iota=1}^{k-1}\mathcal{E}_{\iota}\rt)(s)ds \\
 \lesssim&      \sum_{\iota=0}^{k} \mathcal{E}_{\iota}(0),
\end{split}\ee
since
\bee\label{}\begin{split}
&  \int_0^t    (1+s)^{ -1}  \mathcal{E}_{k} (s)ds
= \int_0^t (1+ s)^{2(k-1)+1}  \int  \bar\rho_0  \lt(\pl_s^k w\rt)^2  dxds \\
& +  \int_0^t (1+ s)^{2k-1}  \int   \bar\rho_0^\ga \lt(\pl_s^k  w_x \rt)^2  dx ds
 +\int_0^t  (1+ s)^ {2k}  \int  \bar\rho_0    \lt(\pl_s^{k+1} w\rt)^2   dx
  \lesssim \sum_{\iota=0}^{k} \mathcal{E}_{\iota}(0)  .
\end{split}\eee
It finally follows from \ef{ai1} and \ef{ai2} that
\bee\label{}\begin{split}
    \mathcal{E}_k(t)   + \int_0^t   \int \lt[(1+ s)^{2k+1}     \bar\rho_0    \lt(\pl_s^{k+1} w\rt)^2  +  (1+ s)^{2k-1}      \bar\rho_0^\ga    \lt(\pl_s^{k } w_x\rt)^2   \rt] dx ds
     \lesssim \sum_{\iota=0}^k \mathcal{ E}_{\iota}(0) .
\end{split}\eee This finishes the proof of Lemma \ref{lem53}. $\Box$

\subsection{Verification of the a priori assumption}
In this subsection, we prove the following lemma.
\begin{lem}\label{lem31} Suppose that $\mathcal{E}(t)$ is finite, then it  holds  that
\be\label{lem31est}\begin{split}
  & \sum_{j=0}^3 (1+  t)^{2j} \lt\|\pl_t^j w(\cdot,t)\rt\|_{L^\iy}^2   +  \sum_{j=0}^1 (1+  t)^{2j} \lt\|\pl_t^j w_x(\cdot,t)\rt\|_{L^\iy}^2
   \\
   & \qquad +
  \sum_{
  i+j\le l,\   2i+j \ge 4 } (1+  t)^{2j}\lt\|  \varsigma^{\f{2i+j-3}{2}}\pl_t^j \pl_x^i w(\cdot,t)\rt\|_{L^\iy}^2  \lesssim \mathcal{E}(t).
\end{split}\ee
\end{lem}
Once this lemma is proved, the a priori assumption \ef{apriori} is then verified and the proof of Theorem \ref{mainthm} is completed, since it follows from  the nonlinear weighted energy estimate \ef{energy}    and  the elliptic estimate \ef{ellipticestimate} that
$$\mathcal{E}(t)\lesssim \mathcal{E}(0), \ \  t\in [0,T].$$

\noindent {\bf Proof}. We first note that $\mathcal{E}_{j,0}\lesssim \mathcal{E}_j$ for $j=0,\cdots, l$. It follows from \ef{hardybdry}  and \ef{varsigma} that
\bee\label{}\begin{split}
&\int    \varsigma^{\aa-1}\lt(\pl_t^j w\rt)^2   dx \lesssim
\int     d^{\aa-1}\lt(\pl_t^j w\rt)^2   dx  \lesssim  \int     d^{\aa+1} \lt[ \lt(\pl_t^j  w_x \rt)^2 + \lt(\pl_t^j w\rt)^2 \rt]  dx  \\
 \lesssim  &\int     \varsigma^{\aa+1} \lt[ \lt(\pl_t^j  w_x \rt)^2 + \lt(\pl_t^j w\rt)^2 \rt]  dx
 \lesssim \int  \lt[   \varsigma^{\aa+1} \lt(\pl_t^j  w_x \rt)^2 +  \varsigma^{\aa}\lt(\pl_t^j w\rt)^2 \rt]  dx \\
\le &   (1+ t)^ {-2j }   \mathcal{ E}_{j}  ,
\end{split}\eee
which implies
\bee\label{}\begin{split}
 \mathcal{ E}_{j,0}= (1+t)^{2j} \int \lt[  \varsigma^{\aa+1}\lt(\pl_t^j w_x\rt)^2 +    \varsigma^{\aa-1}\lt(\pl_t^j w\rt)^2   \rt]dx   \lesssim  \mathcal{ E}_{j}.
\end{split}\eee
So, we have
\be\label{newnorm}\begin{split}
\sum_{j=0}^l \lt(\mathcal{ E}_{j}(t) + \sum_{i=0}^{l-j} \mathcal{ E}_{j, i}(t)  \rt)\lesssim   \mathcal{E}(t).
\end{split}\ee
The following embedding (cf. \cite{adams}):
$H^{1/2+\da}(\mathcal{I})\hookrightarrow L^\iy(\mathcal{I}) $
with the estimate
\be\label{half} \|F\|_{L^\iy(\mathcal{I})} \le C(\delta) \|F\|_{H^{1/2+\da}(\mathcal{I})},\ee
for $\da>0$ will be used in the rest of the proof.

It follows from \ef{wsv} and \ef{varsigma} that for $j\le 5+[\aa]-\aa$,
\be\label{t311}\begin{split}
&\lt\|\pl_t^j w \rt\|_{H^{\frac{5-j+[\aa]-\aa}{2}}}^2
=  \lt\|\pl_t^j w \rt\|_{H^{l-j+1-\frac{l-j+1+\aa}{2}}}^2
\lesssim\lt\|\pl_t^j w \rt\|_{H^{ {l-j+1+\aa} , l-j+1}}^2 \\
= &\sum_{k=0}^{l-j+1} \int    d^{\aa+1+l-j}|\pl_x^k \pl_t^j w|^2dx
\lesssim \sum_{k=0}^{l-j+1} \int   \varsigma^{\aa+1+l-j}|\pl_x^k \pl_t^j w|^2dx
\\
\lesssim & \sum_{k=0}^{l-j+1} \int   \varsigma^{\aa+k}|\pl_x^k \pl_t^j w|^2dx
\le   (1+t)^{-2j}\lt(\mathcal{ E}_{j}(t) + \sum_{k=1}^{l-j } \mathcal{ E}_{j,k}(t)\rt)\le (1+t)^{-2j}  \mathcal{E}(t).
\end{split}\ee
This, together with \ef{half}, gives
\bee\label{}\begin{split}
\sum_{j=0}^3 (1+  t)^{2j} \lt\|\pl_t^j w\rt\|_{L^\iy}^2   +  \sum_{j=0}^1 (1+  t)^{2j} \lt\|\pl_t^j w_x\rt\|_{L^\iy}^2 \lesssim \mathcal{E}(t).
\end{split}\eee

To bound the third term on the left-hand side of \ef{lem31est}, we denote
$$\psi :=  \varsigma^{\f{2i+j-3}{2}}\pl_t^j \pl_x^i w.$$
In the following, we prove that $\lt\|\psi\rt\|_{L^\iy}^2
\lesssim     (1+t)^{-2j}\mathcal{E}(t) $ by separating the cases when $\alpha$ is or is not an integer.

{\bf Case 1} ($\aa\neq [\aa]$). When $\aa$ is not an integer,  we choose  $  \varsigma^{2(l-i-j)+ \aa-[\aa]}$ as the spatial weight. A simple calculation yields
\bee\label{}\begin{split}
 \lt|\pl_x   \psi \rt| \lesssim\lt|  \varsigma^{\f{2i+j-3}{2}}\pl_t^j \pl_x^{i+1} w\rt|
 + \lt|  \varsigma^{\f{2i+j-3}{2}-1}\pl_t^j \pl_x^i w\rt|,
 \end{split}\eee
\bee\label{}\begin{split}
 \lt|\pl_x^2 \psi \rt| \lesssim\lt|  \varsigma^{\f{2i+j-3}{2}}\pl_t^j \pl_x^{i+2} w\rt|
 + \lt|  \varsigma^{\f{2i+j-3}{2}-1}\pl_t^j \pl_x^{i+1} w\rt| + \lt|  \varsigma^{\f{2i+j-3}{2}-2}\pl_t^j \pl_x^i w\rt|,
\end{split}\eee
$$\cdots\cdots$$
\be\label{l313}\begin{split}
 &\lt|\pl_x^k \psi \rt| \lesssim  \sum_{p=0}^k \lt|  \varsigma^{\f{2i+j-3}{2}-p}\pl_t^j \pl_x^{i+k-p} w\rt|   \ \ {\rm for} \ \   k=1, 2, \cdots, l-j+1-i.
\end{split}\ee
It  follows from \ef{l313} that for $1\le k \le l+1-i-j$,
\bee\label{}\begin{split}
& \int   \varsigma^{2(l-i-j)+ \aa-[\aa]  } \lt|\pl_x^k \psi \rt|^2 dx \lesssim   \int  \sum_{p=0}^k    \varsigma^{\aa+l-j+1-2p}\lt|\pl_t^j  \pl_x^{i+k-p} w\rt|^2 dx \\
 \lesssim &   \int   \varsigma^{l-i-j+1-k}  \sum_{p=0}^1     \varsigma^{\aa+i+k-2p} \lt|\pl_t^j  \pl_x^{i+k-p} w\rt|^2 dx + \int  \sum_{p=2}^k    \varsigma^{\aa+l-j+1-2p}\lt|\pl_t^j  \pl_x^{i+k-p} w\rt|^2 dx
 \\
 \lesssim & (1+t)^{-2j} \mathcal{E}_{j,i+k-1}+ \int  \sum_{p=2}^k    \varsigma^{\aa+l-j+1-2p}\lt|\pl_t^j  \pl_x^{i+k-p} w\rt|^2 dx.
\end{split}\eee
To bound the 2nd term on the right-hand side of the inequality above, notice that
\be\label{3.14}\begin{split}
&\aa+l-j+1-2p \\
= &2( l+1-i-j- k)   +2(k-p) + (\aa-[\aa])+(2i+j-3) -2 > -1
\end{split}\ee
for $p\in [2, k]$, due to $\aa\neq [\aa]$ and $2i+j\ge 4$. We then have, with the aid of \ef{hardybdry} and \ef{varsigma}, that for $p\in [2, k]$,
\bee\label{}\begin{split}
  &\int   \varsigma^{\aa+l-j+1-2p}\lt|\pl_t^j  \pl_x^{i+k-p} w\rt|^2 dx
  \lesssim \int     d^{\aa+l-j+1-2p}\lt|\pl_t^j  \pl_x^{i+k-p} w\rt|^2 dx \\
 \lesssim &    \int       d^{\aa+l-j+1-2p+2} \sum_{\iota=0}^1 \lt|\pl_t^j  \pl_x^{i+k-p+\iota} w\rt|^2 dx \lesssim \cdots \cdots \\
 \lesssim &  \int       d^{\aa+l-j+1 } \sum_{\iota=0}^p \lt|\pl_t^j  \pl_x^{i+k-p+\iota} w\rt|^2 dx \lesssim \int      \varsigma^{\aa+l-j+1 } \sum_{\iota=0}^p \lt|\pl_t^j  \pl_x^{i+k-p+\iota} w\rt|^2 dx\\
 =&\int \sum_{\iota=0}^p   \varsigma^{(l+1-i-j-k)+(p-\iota)}     \varsigma^{\aa+i+k-p+\iota }  \lt|\pl_t^j  \pl_x^{i+k-p+\iota} w\rt|^2 dx\\
 \lesssim & \sum_{\iota=0}^p \int   \varsigma^{\aa+i+k-p+\iota }  \lt|\pl_t^j  \pl_x^{i+k-p+\iota} w\rt|^2 dx
 \le   \sum_{\iota=i+k-p}^{i+k-1} (1+t)^{-2j} \mathcal{E}_{j,\iota}.
\end{split}\eee
That yields, for $k=1, 2, \cdots, l-j+1-i$,
\bee\label{}\begin{split}
  \int   \varsigma^{2(l-i-j)+ \aa-[\aa]  } \lt|\pl_x^k \psi \rt|^2 dx
 \lesssim & (1+t)^{-2j} \mathcal{E}_{j,i+k-1}+ \sum_{p=2}^k  \sum_{\iota=i+k-p}^{i+k-1} (1+t)^{-2j} \mathcal{E}_{j,\iota} \\
 \lesssim &(1+t)^{-2j} \sum_{\iota=i}^{i+k-1} \mathcal{E}_{j,\iota}.
\end{split}\eee
Therefore,  it follows from  \ef{varsigma}  and \ef{newnorm} that
\bee\label{}\begin{split}
&\lt\|\psi\rt\|_{H^{2(l-i-j)+\aa-[\aa] , \ l+1-i-j }}^2 = \sum_{k=0}^{l+1-i-j} \int    d^{2(l-i-j)+ \aa-[\aa]  } \lt|\pl_x^k \psi \rt|^2 dx
\\
\lesssim & \sum_{k=0}^{l+1-i-j} \int   \varsigma^{2(l-i-j)+ \aa-[\aa]  } \lt|\pl_x^k \psi \rt|^2dx
 \lesssim    \int   \varsigma^{2(l-i-j)+ \aa-[\aa]  } \lt| \psi \rt|^2 dx+  (1+t)^{-2j} \sum_{\iota=i}^{l-j} \mathcal{E}_{j,\iota} \\
 \lesssim  & (1+t)^{-2j} \sum_{\iota=i}^{l-j} \mathcal{E}_{j,\iota}
 \le  (1+t)^{-2j}\mathcal{E}(t).
\end{split}\eee
When $\aa$ is not an integer, $\aa-[\aa]\in(0,1)$. So, it follows from  \ef{half} and \ef{wsv}  that
\bee\label{}\begin{split}
\lt\|\psi\rt\|_{L^\iy}^2 \lesssim \lt\|\psi\rt\|_{H^{1-\frac{\aa-[\aa]}{2}}}^2
\lesssim \lt\|\psi\rt\|_{H^{2(l-i-j)+\aa-[\aa] , \ l+1-i-j }}^2
\lesssim     (1+t)^{-2j}\mathcal{E}(t).
 \end{split}\eee
{\bf Case 2} ($\aa=[\aa]$). In this case $\aa$ is an integer,  we choose  $  \varsigma^{2(l-i-j)+ 1/2 }$ as the spatial weight. As shown in Case 1, we
have   for $1\le k \le l+1-i-j$,
\bee\label{}\begin{split}
  \int   \varsigma^{2(l-i-j)+ 1/2} \lt|\pl_x^k \psi \rt|^2dx
 \lesssim   (1+t)^{-2j} \mathcal{E}_{j,i+k-1}+ \int  \sum_{p=2}^k    \varsigma^{\aa+l-j+1-2p+ 1/2}\lt|\pl_t^j  \pl_x^{i+k-p} w\rt|^2 dx
\end{split}\eee
and
\bee\label{}\begin{split}
 \aa+l-j+1-2p + \frac{1}{2}
=  2( l+1-i-j- k)   +2(k-p)  +(2i+j-3) - \frac{3}{2} \ge  -\frac{1}{2}.
\end{split}\eee
We can use the Hardy inequality to obtain
\bee\label{}\begin{split}
  \int   \varsigma^{2(l-i-j)+ 1/2 } \lt|\pl_x^k \psi \rt|^2 dx
 \lesssim  (1+t)^{-2j} \sum_{\iota=i}^{i+k-1} \mathcal{E}_{j,\iota}, \ \ k=1, 2, \cdots, l-j+1-i
\end{split}\eee
and
\bee\label{}\begin{split}
&\lt\|\psi\rt\|_{H^{2(l-i-j)+1/2 , \ l+1-i-j }}^2
 \lesssim    (1+t)^{-2j} \sum_{\iota=i}^{l-j} \mathcal{E}_{j,\iota}.
\end{split}\eee
Therefore, it follows from  \ef{half} and \ef{wsv}  that
\bee\label{}\begin{split}
\lt\|\psi\rt\|_{L^\iy}^2 \lesssim \lt\|\psi\rt\|_{H^{3/4}}^2
\lesssim \lt\|\psi\rt\|_{H^{2(l-i-j)+1/2 , \ l+1-i-j }}^2
\lesssim     (1+t)^{-2j}\mathcal{E}(t).
 \end{split}\eee
This completes the proof of  Lemma \ref{lem31}. $\Box$

\section{Proof of Theorem \ref{mainthm1}}
In this section, we prove Theorem \ref{mainthm1}.
First, it follows from \ef{ldv}, \ef{ld}, \ef{bar} and \ef{212} that for $(x,t)\in \mathcal{I}\times[0,\iy)$,
$$
 \rho(\eta(x, t), t)-\bar\rho(\bar\eta(x, t), t)
 =\frac{\bar\rho_0(x)}{\eta_x(x, t)}-\frac{\bar\rho_0(x)}{\bar\eta_x(x, t)}
 =-\bar\rho_0(x)\frac{w_x(x, t)+h(t)}{(\tilde \eta_x+w_x)\bar\eta_x(x, t)}
$$
and
$$
u(\eta(x, t), t)-\bar u (\bar\eta(x, t), t)=w_t(x, t)+xh_t(t) . $$
Hence, by virtue of \ef{basic}, \ef{212}, \ef{2} and \ef{h}, we have, for $x \in \mathcal{I}$ and $t\ge 0$,
$$
|\rho(\eta(x, t), t)-\bar\rho(\bar\eta(x, t), t)|\lesssim \lt(A-Bx^2\rt)^{\frac{1}{\ga-1}}(1+t)^{-\frac{2}{\ga+1}}\lt(\sqrt{\mathcal{E}(0)}+
(1+t)^{-\frac{\ga}{\ga+1}}\ln(1+t)\rt) $$
and
$$|u(\eta(x, t), t)-\bar u (\bar\eta(x, t), t)|\lesssim (1+t)^{-1}\sqrt{\mathcal{E}(0)}+(1+t)^{-\frac{2\ga+1}{\ga+1}} \ln(1+t).$$
Then \ef{1'} and \ef{2'} follow.  It follows from \ef{vbs} and \ef{212} that
\bee\begin{split}
x_+(t) = &\eta\lt(\bar x_+(0),t\rt)=\lt(\tilde{\eta}+ w\rt)\lt(\bar x_+(0),t\rt)=\lt(\bar{\eta}+ xh + w\rt)\lt(\bar x_+(0),t\rt) \\
 =& \sqrt{ {A}/{B}} \lt((1+t)^{\frac{1}{\ga+1}}+h(t)\rt)+w\lt(\sqrt{ {A}/{B}}, t\rt), \end{split}\eee
which, together with  \ef{h} and \ef{2}, implies  that  for $t\ge 0$,
\bee\label{}
\sqrt{ {A}/{B}}(1+t)^{\frac{1}{\ga+1}}-C\sqrt{\mathcal{E}(0)}\le x_+(t)
\le \sqrt{ {A}/{B}}\lt((1+t)^{\frac{1}{\ga+1}}+C(1+t)^{-\frac{1}{\ga+1}}\rt)+C\sqrt{\mathcal{E}(0)}. \eee
Similarly, we have  for $t\ge 0$,
\bee\label{}
- \sqrt{ {A}/{B}}\lt((1+t)^{\frac{1}{\ga+1}}+C(1+t)^{-\frac{1}{\ga+1}}\rt)-C\sqrt{\mathcal{E}(0)}\le x_-(t)\le  - \sqrt{ {A}/{B}}(1+t)^{\frac{1}{\ga+1}}+ C\sqrt{\mathcal{E}(0)}. \eee
Thus, \ef{3'} follows from  the smallness of $\mathcal{E}(0)$. Notice that
for $k=1, 2, 3$,
$$\frac{d^k x_{\pm}(t)}{dt^k}=\pl_t^ k\tilde \eta \lt(\pm \sqrt{ {A}/{B}}, t\rt)+ \pl_t^ k w \lt(\pm \sqrt{ {A}/{B}}, t\rt).$$
Therefore, $\ef{decay}$ and \ef{2} implies \ef{4'}. $\Box$

\section{Appendix}
 In this appendix, we prove \ef{decay} and \ef{h}. We may write \ef{mt}  as the following system
\be\label{ansatzeq}\lt\{\begin{split}
 &h_t= z,\\
 &z_t= -z -   \lt[ { \bar\eta_x ^{-\ga}}  -   {(\bar\eta_x+h)^{- \ga}}\rt] /({\ga+1}) - \bar\eta_{xtt},\\
 &(h,z)(t=0)=(0,0).
\end{split}\rt.\ee
Recalling that $\bar\eta_x (t)=(1+t)^{\frac{1}{\ga+1}}$, thus  $\bar\eta_{xtt}<0$. A simple phase plane analysis shows that there exist $0<t_0<t_1<t_2$ such that, starting from $(h, z)=(0, 0)$ at $t=0$, $h$ and $z$ increases in the interval $[0, t_0]$ and $z$ reaches its positive maxima at $t_0$; in the interval $[t_0, t_1]$, $h$ keeps increasing and reaches its maxima at $t_1$,  $z$ decreases from its positive maxima to 0; in the interval $[t_1, t_2]$,
both $h$ and $z$ decrease, and $z$ reaches its negative minima at $t_2$; in the interval $[t_2, \infty)$, $h$ decreases and $z$ increases, and $(h, z)\to (0, 0)$ as $t\to \infty$. This can be summarized as follows: \\
$$z(t) \uparrow_0 , \  \ h(t) \uparrow_0  ,  \ \ t\in [0,t_0]$$  $$ z(t) \downarrow_0,  \ \  h(t)  \uparrow,  \ \  t\in[t_0,t_1] $$
$$ z(t) \downarrow^0,  \ \  h(t)  \downarrow,  \ \  t\in[t_1,t_2] $$
$$ z(t)  \uparrow^0,  \ \  h(t)  \downarrow_0,  \ \  t\in[t_2,\iy). $$
It follows from the above analysis that there exists a  constant $C=C(\ga, M)$ such that
\be\label{boundforh}
0\le h(t) \le C  \ \ {\rm for} \ \  t\ge 0.\ee
In view of \ef{ansatz}, we then see that  for some constant $K>0$,
$$\lt(1 +   t \rt)^{{1}/({\ga+1})} \le \tilde \eta_{x} \le K \lt(1 +   t \rt)^{{1}/({\ga+1})}.$$

To derive the decay property,
we may rewrite \ef{mt} as
\be\label{gt}\begin{split}
&\tilde \eta_{x tt}+\tilde \eta_{xt}- \tilde \eta_x^{-\ga}/(\ga+1)=0,\\
&\tilde \eta_x (t=0)= 1, \ \  \tilde \eta_{xt}(t=0)={1 }/({\ga+1}) .
\end{split}\ee
  Then, we have by solving \ef{gt} that
\be\label{k0}\begin{split}
\tilde \eta_{xt}(t)=\f{ 1 }{\ga+1}e^{-t}  + \f{ 1 }{\ga+1} \int_0^t e^{-(t-s)}\tilde \eta_x^{-\ga}(s) ds \ge 0.
\end{split}\ee
Next, we use the mathematical induction to prove $\ef{decay}_{2}$.
First, it follows from \ef{k0} that
\bee\begin{split}
({\ga+1})   \tilde \eta_{xt}(t)  = &e^{-t}   + \int_0^{t/2} e^{-(t-s)}\tilde \eta_{x}^{-\ga}(s) ds +   \int_{t/2}^t e^{-(t-s)}\tilde \eta_{x}^{-\ga}(s) ds \\
\le & e^{-t} +  e^{-t/2}\int_0^{t/2}(1+  s)^{-\f{\ga}{\ga+1}}ds
+\lt(1+  {t}/{2} \rt)^{-\f{\ga}{\ga+1}} \int_{t/2}^t e^{-(t-s)}  ds \\
\le & e^{-t} + \f{ e^{-t/2} }{1+\ga} \lt(1+  {t}/{2}\rt)^{{1}/({\ga+1})}
+\lt(1+ {t}/{2} \rt)^{-{\ga}/({\ga+1})} \\
\le  &  C \lt(1+ {t} \rt)^{-{\ga}/({\ga+1})}, \ t\ge 0,
\end{split}\eee
for some positive constant $C$ independent of $t$. This proves $\ef{decay}_{2}$ for $k=1$.
For  $2\le m \le n  $ for a fixed positive integer $n$, we make the induction hypothesis that $\ef{decay}_{2}$ holds for all $k=1,2,\cdots,m-1$, that is,
\be\label{induction} \lt|\f{d^k\tilde \eta_{x}(t)}{dt^k}\rt| \le C(m)\lt(1 +   t \rt)^{\frac{1}{\ga+1}-k},   \ \ k= 1, 2, \cdots, m-1.\ee
It suffices to prove $\ef{decay}_{2}$ holds  for $k=m$.
We derive from \ef{gt} that
$$
\f{d^{m+1}\tilde \eta_{x}}{dt^{m+1}}(t)+\f{d^{m}\tilde \eta_{x}}{dt^{m}}(t)- \frac{1}{\ga+1}\frac{d^{m-1}\tilde \eta_x^{-\ga}}{dt^{m-1}}(t)=0 , $$
so that
\be\label{gttt}
   \f{d^{m}\tilde \eta_{x}}{dt^{m}}(t) = e^{-t} \f{d^{m}\tilde \eta_{x}}{dt^{m}}(0)  + \frac{1}{\ga+1} \int_0^{t} e^{-(t-s)}\frac{d^{m-1}\tilde \eta_x^{-\ga}}{ds^{m-1}}(s)ds, \ee
where $\f{d^{m}\tilde \eta_{x}}{dt^{m}}(0) $ can be determined by the equation inductively. To bound the last term on the right-hand side of \ef{gttt}, we are to derive that
\be\label{estimate1}
 |\pl_t^r(\tilde \eta_x^{-1})|( t)\le C(m) (1+t)^{-\frac{1}{\ga+1}-r},  \  \  0\le r\le m-1, \ee
for some constant $C(m)$ depending only on $\ga$, $M$ and $m$. First, \ef {estimate1} is true for $r=0$ in view of $(2.12)_1$. For $1 \le r\le m-1$, note that
\begin{align*}\label{}
 \pl_t^r(\tilde \eta_x^{-1})&=\pl_t^{r-1}(\tilde \eta_x^{-2}\tilde \eta_{xt})
=\sum_{i=0}^{r-1} C_{r-1}^i \pl_t^{i}(\tilde \eta_x^{-2}) \pl_t^{r-i}(\tilde \eta_x)\notag\\ & = \sum_{i=0}^{r-1} C_{r-1}^i \pl_t^{i}(\tilde \eta_x^{-2}) \pl_t^{r-i}(\tilde \eta_x)\notag\\
&=\sum_{i=0}^{r-1} C_{r-1}^i  \lt(\sum_{j=0}^{i} C_i^j \pl_t^{j}(\tilde \eta_x^{-1}) \pl_t^{i-j}(\tilde \eta_x^{-1})\right)\pl_t^{r-i}(\tilde \eta_x).
\end{align*}
Then, \ef{estimate1} can be  proved by an iteration, with the aid of \ef{induction}.
Notice that
\bee\begin{split}\label{xx1}
 \pl_t^{m-1}( \tilde \eta_{x}^{-\ga})&=-\ga\pl_t^{m-2}\lt( \tilde \eta_{x}^{-(\ga+1)}\tilde \eta_{xt}\rt)=-\ga\sum_{i=0}^{m-2}C_{m-2}^i\pl_t^{i}\lt( \tilde \eta_{x}^{-(\ga+1)}\rt)\lt(\pl_t^{m-1-i}\tilde \eta_{xt}\rt) \notag\\
&=\ga(\ga+1)\sum_{i=0}^{m-2}C_{m-2}^i\lt[\sum_{j=0}^iC_i^j\pl_t^{j}\lt( \tilde \eta_{x}^{-\ga}\rt)\pl_t^{i-j}\lt(\tilde \eta_x^{-1}\rt)\rt]\lt(\pl_t^{m-1-i}\tilde \eta_{xt}\rt). \end{split}\eee
It therefore follows from  \ef{estimate1} and  \ef{induction}  that
\be\label{xx2}
|\pl_t^{m-1}(\tilde\eta_x^{-\ga})|\le C_1(m)(1+t)^{-\frac{\ga}{\ga+1}-(m-1)}\ee
for some constant $C_1(m)$ independent of $t$. This, together with \ef{gttt}, proves that $\ef{decay}_{2}$  is also true for $k=m$, and completes the proof of $\ef{decay}_{2}$.

Finally, we prove  the decay estimate for $h$. We may write the equation for $h$ as
\be\label{zuihou}
h_t+\frac{1}{\ga+1}(1+t)^{-\frac{\ga}{\ga+1}}\lt[1-\lt(1+ h(1+t)^{-\frac{1}{\ga+1}}\rt)^{-\ga}\rt]=-\tilde \eta_{xtt}.
\ee
Notice that
$$\lt(1+ h(1+t)^{-\frac{1}{\ga+1}}\rt)^{-\ga}\le 1-\ga h (1+t)^{-\frac{1}{\ga+1}}+\frac{\ga(\ga+1)}{2} h^2 (1+t)^{-\frac{2}{\ga+1}},$$
due to $h\ge 0$. We then obtain, in view of $\ef{decay}_{2}$, that
$$
h_t+\frac{\ga}{\ga+1}(1+t)^{-1}h\le \frac{\ga}{2}(1+t)^{-\frac{\ga+2}{\ga+1}}h^2+C (1+t)^{\frac{1}{\ga+1}-2}.$$
Thus,
\be\label{7201}
h(t)\le C(1+t)^{-\frac{\ga}{\ga+1}} \int_0^t \lt((1+s)^{-\frac{2}{\ga+1}}h^2(s)+(1+s)^{-1}\rt)ds .\ee
We use an iteration to prove \ef{h}.
First, since $h$ is bounded due to \ef{boundforh}, we have
\be\label{7202}h(t)\le C(1+t)^{-\frac{\ga}{\ga+1}} \int_0^t (1+s)^{-\frac{2}{\ga+1}}ds \le C(1+t)^{-\frac{1}{\ga+1}}. \ee
Substituting this into \ef{7201}, we obtain
$$h(t) \le C(1+t)^{-\frac{\ga}{\ga+1}} \int_0^t \lt((1+s)^{-\frac{4}{\ga+1}}+(1+s)^{-1}\rt)ds ;$$
which implies
\bee\label{}h(t) \le \lt\{ \begin{split} & C (1+t)^{-\frac{\ga}{\ga+1}}\ln(1+t) & {\rm if} \ \  \ga \le 3,\\
& C (1+t)^{-\frac{3}{\ga+1}}  & {\rm if} \ \  \ga > 3. \end{split}\rt. \eee
If $\ga\le 3$, then the first part of \ef{h} has been proved. If $\ga>3$, we repeat this procedure and obtain
\bee\label{}h(t) \le \lt\{ \begin{split} & C (1+t)^{-\frac{\ga}{\ga+1}}\ln(1+t) & {\rm if} \ \  \ga \le 7,\\
& C (1+t)^{-\frac{7}{\ga+1}}  & {\rm if} \ \  \ga > 7. \end{split}\rt. \eee
For general $\gamma$, we repeat this procedure $k$ times with $k= \lceil\log_2(\ga+1)\rceil  $ satisfying
$\sum_{j=0}^{k} 2^j\ge \ga$ to  obtain
$$h(t) \le C (1+t)^{-\frac{\ga}{\ga+1}}\ln(1+t). $$
This, together with \ef{boundforh}, proves the first part of \ef{h}, which in turn  implies the second part of \ef{h}, by virtue of  \ef{zuihou} and $\ef{decay}_2$.

\vskip 0.5cm
\centerline{\bf Acknowledgement}

Luo's research was  supported in part by NSF under grant DMS-1408839,  Zeng's research was supported in part by NSFC under grant \#11301293/A010801.

\bibliographystyle{plain}

\noindent {Tao Luo}\\
Department of Mathematics and Statistics\\
Georgetown University,\\
Washington, DC, 20057, USA. \\
Email: tl48@georgetown.edu\\

\noindent Huihui Zeng\\
Mathematical Sciences Center\\
Tsinghua University\\
Beijing, 100084, China.\\
E-mail: hhzeng@mail.tsinghua.edu.cn

\end{document}